\documentclass[11pt,reqno]{amsart}
\usepackage[table,xcdraw]{xcolor}
\usepackage[colorlinks=true,linkcolor=black,urlcolor=black,citecolor=black]{hyperref}
\usepackage{a4wide}
\usepackage{euscript,amsmath,amssymb,amsbsy,mathabx,esint}
\usepackage[arrow,matrix,curve]{xy}\usepackage{tikz}\usetikzlibrary{patterns}
\usepackage{array,longtable,multirow,diagbox,graphicx,enumitem,url,hhline,relsize,scalefnt}
\usepackage{float,caption,capt-of}
\newcounter{msct}[section]\renewcommand{\themsct}{\thesection.\arabic{msct}}
\newenvironment{m-theorem}{\vskip3pt\refstepcounter{msct}\trivlist \itemindent 0pt%
\item[\hskip\labelsep\bf Theorem \themsct]\it\ignorespaces}{\endtrivlist\vskip2pt}
\newenvironment{m-proposition}{\vskip3pt\refstepcounter{msct}\trivlist \itemindent0pt%
\item[\hskip\labelsep\bf Proposition \themsct]\it\ignorespaces}{\endtrivlist\vskip2pt}
\newenvironment{m-corollary}{\vskip3pt\refstepcounter{msct}\trivlist \itemindent 0pt%
\item[\hskip\labelsep\bf Corollary \themsct]\it\ignorespaces}{\endtrivlist\vskip2pt}
\newenvironment{m-lemma}{\vskip3pt\refstepcounter{msct}\trivlist \itemindent 0pt%
\item[\hskip\labelsep\bf Lemma \themsct]\it\ignorespaces}{\endtrivlist\vskip2pt}
\newenvironment{m-definition}{\vskip3pt\refstepcounter{msct}\trivlist \itemindent0pt%
\item[\hskip\labelsep\bf Definition \themsct]\ignorespaces}{\endtrivlist\vskip2pt}
\newenvironment{m-notation}{\vskip3pt\refstepcounter{msct}\trivlist \itemindent0pt%
\item[\hskip\labelsep\bf Notation \themsct]\ignorespaces}{\endtrivlist\vskip3pt}
\newenvironment{m-example}{\vskip3pt\refstepcounter{msct}\trivlist \itemindent0pt%
\item[\hskip\labelsep\bf Example \themsct]\ignorespaces}{\endtrivlist\vskip3pt}
\newenvironment{m-remark}{\vskip3pt\refstepcounter{msct}\trivlist \itemindent0pt%
\item[\hskip\labelsep\bf Remark \themsct]\ignorespaces}{\endtrivlist\vskip3pt}
\newenvironment{m-question}{\vskip3pt\refstepcounter{msct}\trivlist \itemindent0pt%
\item[\hskip\labelsep\bf Issue.]\ignorespaces}{\endtrivlist\vskip3pt}
\newenvironment{thm-nono}[1]{\vskip3pt\trivlist \itemindent 0pt %
\item[\hskip\labelsep{\bf Theorem}~#1]\it\ignorespaces}{\endtrivlist\vskip3pt}
\newenvironment{lm-nono}[1]{\vskip3pt\trivlist \itemindent0pt%
\item[\hskip\labelsep{\bf Lemma}~#1]\it\ignorespaces}{\endtrivlist\vskip3pt}
\newenvironment{conj-nono}[1]{\vskip3pt\trivlist \itemindent0pt%
\item[\hskip\labelsep{\bf Conjecture}~#1]\it\ignorespaces}{\endtrivlist\vskip3pt}
\newenvironment{m-thank}{\vskip3pt\trivlist \itemindent0pt%
\item[\hskip\labelsep\it Acknowledgments]\ignorespaces}{\endtrivlist\vskip3pt}
\newenvironment{m-proof}{\vskip2pt\trivlist \itemindent0pt%
\item[\hskip\labelsep\it Proof.]\ignorespaces}{\hfill$\Box$\endtrivlist\vskip3pt}%
\newenvironment{m-asmp}{\vskip3pt\trivlist \itemindent0pt%
\item[\hskip\labelsep\bf Assumption.]\ignorespaces}{\hfill\endtrivlist\vskip3pt}%

\newcounter{meqn}[section]\renewcommand{\themeqn}{\thesection.\arabic{meqn}}
\newenvironment{m-eqn}[1]{\vskip5pt\refstepcounter{meqn}\trivlist\itemindent0pt%
\item[]\ignorespaces\hfill$\displaystyle #1$\hfill\hbox{\rm(\themeqn)}}{\endtrivlist\vskip3pt}

\newcommand{\bibauth}[2]{\textrm{{#1}~{#2}},}
\newcommand{\bibtitl}[1]{\textit{#1},}
\newcommand{\bibjnyp}[4]{\textrm{#1} \textbf{#2} (#3), #4.}

\newcommand{\bibbook}[4]{\textit{#1}. {#2} {#3}, {#4}.}

\numberwithin{equation}{section}

\let\mt\mapsto
\let\lmt\longmapsto

\font\tenmsa=msam10 %
\newcommand\hdashpiece{%
{\vrule height2.75pt depth-2.35pt width2.3pt \kern1.7pt}}%
\newcommand\hdashpieces{%
{\hdashpiece\hdashpiece\hdashpiece\hdashpiece}}%
\newcommand\dashto{\mathrel{%
\hdashpiece\hdashpiece\kern-0.4pt\hbox{\tenmsa K}}}%
\newcommand\dashar{\mathrel{%
\hdashpieces\kern-0.4pt\hbox{\tenmsa K}}}%

\let\euf\EuScript 

\let\mbb\mathbb

\let\bsymb\boldsymbol 
\DeclareFontFamily{OT1}{rsfs}{}
\DeclareFontShape{OT1}{rsfs}{n}{it}{<->rsfs10}{}
\DeclareMathAlphabet{\crl}{OT1}{rsfs}{n}{it}

\let\aph\alpha 
\let\blt\bullet

\let\disp\displaystyle
\let\dta\delta

\let\gma\gamma
\let\ges\geqslant   \let\les\leqslant
\let\kp\kappa
\let\lda\lambda

\let\nit\noindent

\let\unbar\underbar

\let\tld\tilde

\let\srel\stackrel

\let\veps\varepsilon

\newcommand{\aaa}{a}
\newcommand{\bb}{b}

\newcommand\bone{{1\kern-0.57ex\rm l}}
\newcommand{\bvp}{{\rm BVP}}\newcommand{\ivp}{{\rm IVP}}

\newcommand{\CC}{C}

\newcommand{\exc}{{\rm ex}}

\newcommand{\infl}{{i\!n\!f\!l}}
\newcommand{\mdl}{{m\!i\!d}}

\newcommand{\ort}{\mathrel{{\vrule width5pt height0.5pt depth0pt\vrule width0.5pt height7pt depth0pt\,}}}
\newcommand{\rav}{\mathrel{{\vrule width5pt height6pt depth-5.5pt\vrule width0.55pt height6pt depth0pt\,}}}
\newcommand{\ouset}[3]{\underset{#1}{\overset{#2}{#3}}}
\newcommand{\rd}{{\rm d}}
\newcommand{\mm}{m}
\newcommand{\mn}{{\rm min}}
\newcommand{\nn}{n}
\newcommand{\pp}{p}
\newcommand{\rr}{r}

\newcommand{\turn}{{t\!u\!r\!n}}
\newcommand{\up}{u_+}
\newcommand{\UP}{U_+}\newcommand{\UM}{U_-}
\newcommand{\uu}{u}\newcommand{\UU}{U}
\newcommand{\vu}[1]{u_{#1}}
\newcommand{\ww}{w}
\newcommand{\yp}{y_+}\newcommand{\ym}{y_-}
\newcommand{\yy}{y}\newcommand{\YY}{Y}
\newcommand{\zz}{z}


\definecolor{lgray}{HTML}{e0eeee}

\begin{document}

\title[Mixed approach to the Emden-Fowler equation]{An analytical-numerical approach\\ to the Emden-Fowler equation}

\author{Mihai Halic}

\begin{abstract}
We investigate the EF-equation $\;\uu''=-x^r\uu^\pp,\,\uu(0)=1,\,\uu'(0)=0$, with roots in astrophysics, and study the qualitative and quantitative dependence of its solution on the parameters; the analytical work is paralleled by numerical simulations. A special attention is given to estimating the first zero $x_0$ of $\uu$ in terms of $\rr, \pp$. The results are used to address two apparently new issues: first, we solve EF backwards starting from $x_0$; second, we transform it into an overdetermined boundary value problem and decide when is this solvable.
\end{abstract}

\maketitle

\section*{Introduction}

The initial value problem (\ivp), 
\begin{m-eqn}{
\uu''=- x^r\cdot\uu^\pp,\quad\uu(0)=1,\;\uu'(0)=0, 
}\label{eq:ivp}
\end{m-eqn}
is known as the Emden-Fowler equation and has relevance in cosmology and astrophysics~\cite{emden,cha}.  Fowler~\cite{fow1,fow2} was among the first to study the analytical properties of the solution. The polytropic index $\pp$ typically ranges from $1$ to $5$; values greater than $5$ occur for modelling isothermal gaseous spheres~\cite{tar+sak} and early universe scenarios. 

There is a massive amount of literature dedicated to the generalized Emden-Fowler-Lane equation, obtained by replacing $x^r$ with an arbitrary function, it is basically impossible to cite all related references. A survey reveals that they are divided into theoretical and numerical approaches, each being technical in its own way. On the theoretical side, we mention~\cite{ram,tal,tal2}; a lengthy, comprehensive review of the analytical results for $\,U''=-a(x)U^\pp,$ with $a(x)$ non-negative, can be found in~\cite{wong}. Concerning numerical solutions, see~\cite{ramos,sita,wrd,wwg} for power series approaches, \cite{he,liao,sws} for homotopy, and~\cite{neural} for neural network methods. 

Mixed \emph{analytical and computational} techniques, where one uses in parallel analytic methods to understand the properties of the solution and computer software to derive numerical information based on the analytical outcome, seem to be rare. Concretely, the solution of~\eqref{eq:ivp} is a differentiable function and one is interested in geometric information about its graph: shape, slope, boundary layer phenomena, $X$-intercept, dependence on parameters. This kind of approach is challenging because it should be at the same time \emph{analytically rigorous and simple enough} to allow extracting explicit, correct estimates. 
In this article, we take a step in this direction, for $\rr>0, \pp\ges1$. 

On the analytical side, we use upper/lower envelope techniques to clarify the dependence of the solution $\uu_\exc$ on parameters. On the computational side, we confront the theoretical results against numerical data, probe their validity, and sharpen them. 
Thus we focus on the relationship (geometry of $\uu_\exc$)---(parameters $\rr, \pp$), which is not covered by numerical methods, rather than on high precision computations, for which there are accurate algorithms. 

The main novelty in our approach is the application of Newton's method to approximate $\uu_\exc$ by a straight line, from the `turning point' (where $\uu_\exc'''=0$) down to its (first) zero $x_0(\rr,\pp)$. It sharply contrasts the frequently used power series approximations. The simplicity of this construction allows estimating the value of $x_0(\rr,\pp)$, relevant to physics, in terms of $\rr, \pp$. Such estimate seems to be missing. 
We verify the correctness of our procedure in two different ways. First, we solve~\eqref{eq:ivp} backwards, starting from $x_0(\rr,\pp)$, to check at what extent we recover the initial values at $x=0$. Second, we address the following \eqref{eq:ef}-related \bvp. 
\\[1ex]\centerline{
\begin{minipage}{0.85\textwidth}\raggedright{
Given $\zz>1$, determine the relationship between the \emph{unknown} parameters $\rr, \pp>2$, such that the overdetermined equation below is solvable: 
\vskip.5ex\centerline{$\uu''=-x^\rr\cdot\uu^\pp,\;\uu(0)=1,\,\uu(\zz)=0,\,\uu'(0)=0.$}
}\end{minipage}
}\\[1ex]  
\nit We are aware of \emph{no references} dealing with this issue. Nevertheless, it's easy to imagine a (physical) scenario leading to it: one is able to compute or measure the zero of the (physically relevant) function $\uu$ (e.g. the radius $\zz$ of a gaseous star) but doesn't know the parameters (inside the star) which determine it. (Thus, the excessive condition $\uu'(0)=0$ is imposed by physical considerations.) The answer turns out to be quite simple (cf.~\S\ref{ssct:zx0}): 
\\[.5ex]\centerline{$\disp 
\pp\approx2\cdot\bigg[\frac{{\zz}^{\frac{\rr+2}{\rr+1}}}{{[(\rr+1)(\rr+2)]}^{\frac{1}{\rr+1}}}-1\bigg]\cdot{(\rr+1)}^{\frac{\rr+3}{\rr+2}}.$
}\\[.5ex] 
We verify that this relationship is in agreement with numerical tests, for values $\zz=1.001\,\text{--}\,500$.



\section{Analytical preparations}\label{sct:theory}

We wish to determine approximate solutions to non-integrable ODEs $\yy''=F(x,\yy)$ by replacing $F$ with a function $G$ such that: 
\begin{itemize}[leftmargin=3ex]
\item $|G-F|$ is small. (So the solution of the new equation approximates the initial one.)
\item The equation $\yy''=G(x,\yy)$ is explicitly integrable.
\end{itemize}
The first requirement is loose, but the second is restrictive. Actually, to estimate the accuracy of the approximate solutions, we need two functions $G_\pm$ satisfying the conditions above, such that the corresponding solutions $\yy_\pm$ `squeeze' in between the exact solution: $\;\ym\les\yy_\exc\les \yp.$ 
We say that $\yy_\pm$ are \emph{upper/lower envelopes} of the exact solution $\yy_\exc$. The next proposition is the toolbox for checking this matter. 

\begin{m-notation}
Let $(y_0,y_1), (z_0,z_1)$ be pairs of real numbers. By $(y_0,y_1)<(z_0,z_1)$, we mean that $y_j\les z_j, j=0,1$, and at least one inequality is strict. We use similar notation for triples.
\end{m-notation}

\begin{m-proposition}\label{prop:+-} 
\nit{\rm(i)} 
Let $F,G$ be differentiable on $[\xi_0,\xi_1]\times[\Xi_0,\Xi_1]$ and let $\yy, \zz$ be respectively the solutions of: 
\vskip.5ex\hskip-3ex\renewcommand{\arraystretch}{1.15}\textscale{.9}{
\begin{tabular}{clll}
\cline{2-4}
either:&$y'=F(x,y),$&$z'=G(x,z),$& 
{$\,(y_0,F(\xi_0,y_0))<(z_0,G(\xi_0,z_0));$}
\\[0ex] 
{\ivp1}&$y(\xi_0)=y_0,$&$z(\xi_0)=z_0,$&
\\[0ex] \cline{2-4}
or:&$y''=F(x,y),$&$z''=G(x,z),$& \multirow{2}{*}{$(y_0,y_0',F(\xi_0,y_0))<(z_0,z_0',G(\xi_0,z_0))$.}
\\ 
{\ivp2}& $y(\xi_0)=y_0,\,y'(\xi_0)=y_0',$ & $z(\xi_0)=z_0,\,z'(\xi_0)=z_0',$ & 
\\\cline{2-4}
\end{tabular}
}\renewcommand{\arraystretch}{1}\vskip.5ex 
\nit Assume that, for some $\Xi\in(z_0,\Xi_1)$, the following conditions hold true: 
\vskip.5ex\centerline{
{\rm(a)} $F(x,r)\les G(x,r)$, for $r\in[\Xi_0,\Xi]$;\quad 
{\rm(b)} either $F$ or $G$ is non-decreasing in $r$.}\vskip.5ex 

\nit Then we have $y(x)\les z(x)$ on $[\xi_0,\xi]$, where 
$\xi:=\max\{x\mid y(\tau), z(\tau)\in[\Xi_0,\Xi],\,\forall\,\tau\in[\xi_0,x] \}.$ 

\medskip\nit{\rm(ii)} Suppose $(y_0,y_0',F(\xi_0,y_0))=(z_0,z_0',G(\xi_0,z_0))$ and there is $e>0$, such that {\ivp2} above has the following property. For all $0<\veps<e$, the solutions $y_{(\veps)}, z$ with the initial conditions 
$$
(y_{(\veps)}(\xi_0),{y_{(\veps)}}'(\xi_0))=(y_0-\veps,y'_0)
\;\;\text{and}\;\;
(z(\xi_0),z'(\xi_0))=(y_0,y'_0), 
$$ 
satisfy $y_{(\veps)}(x)\les z(x)$ for $x\in[\xi_0,\xi]$, where $\xi$ is \emph{independent} of $\veps$. Then one has 
\vskip.5ex\centerline{$\lim_{\veps\to 0}y_{(\veps)}=y_{(0)}(x)\les z(x),\;\text{for}\;x\in[\xi_0,\xi].$}  
\end{m-proposition}

\begin{m-proof}
\nit(i) Let $\dta(x):=z(x)-y(x)$. The hypothesis $(y_0,F(\xi_0,y_0))<(z_0,G(\xi_0,z_0))$ implies that $\dta(x)>0$, for all $x>\xi_0$ sufficiently close to $\xi_0$. If $\{x\mid \dta(x)<0\}$ is non-empty, it contains an interval $(c,d)\subset(\xi_0,\xi_1)$ with $c$ nearest to $\xi_0$, so $\dta(x)>0$ for $x\in(\xi_0,c)$ and $\dta(c)=0$. Thus, on the one hand, the intermediate value theorem implies that there is $\tau\in(\xi_0,c)$, such that $\dta'(\tau)<0$. On the other hand, assuming $G$ is non-decreasing, we obtain the contradiction: 
\\[.5ex]{$z'(\tau)=G(\tau,z(\tau))\ges G(\tau,y(\tau))\ges F(\tau,y(\tau))=y'(\tau).$} 
(A similar argument works for $F$.)

For the second order \ivp, the argument is analogous. The hypothesis $(y_0,y_1',F(\xi_0,y_0))<(z_0,z_1',G(\xi_0,z_0))$ implies that $\dta(\xi_0), \dta'(\xi_0), \dta''(\xi_0)\ges 0$ and at least one is strictly positive. Thus we have $\dta(x)>0$, for $|x-\xi_0|$ sufficiently small. The proof continues as above. 

\nit(ii) The statement follows from the continuous dependence of the solutions on parameters. 
\end{m-proof}


\subsection{Boundary layer}\label{ssct:bd-layer} 

The terminology is common for two-point {\bvp}s. It reflects that a function $\yy=\yy(x)$ varies steeply near a line $x=\xi$, as near a vertical asymptote. 

\begin{m-definition}\label{def:bl}
Let $I=[\xi_0,\xi_1]$ and $\yy:I\to\mbb R$ be a continuous function. For a closed sub-interval $J\subset I$, we denote $|J|$ its length and ${\rm Var}_J(\yy):=\max_{x,x'\in J}|\yy(x)-\yy(x')|$ the total variation. We say that $\yy$ has \emph{boundary layer} at $\xi\in I$ in the following situation: 
\begin{itemize}[leftmargin=3ex]
\item 
There is $\veps>0$, such that $\yy$ is monotone on $J=[\xi-\veps,\xi]\subset I$ (or $J=[\xi,\xi+\veps]\subset I$); 
\item 
${\rm Var}_J(\yy)>50\cdot|J|$ and ${\rm Var}_J(\yy)>0.5\cdot{\rm Var}_{[\xi_0,\xi]}(\yy)$. (One may replace $0.5$ as is convenient.)
\end{itemize}
\end{m-definition}
\nit Suppose $\yy$ is differentiable on $J$. The first condition is satisfied if $|\yy'(x)|>50,\;\forall x\in J$; if $\yy$ is also convex or concave, it's enough $|\yy'(\xi\mp\veps)|>50$. 
This  explains the term `boundary' layer: since $\tan(89^\circ)\approx50$, the condition $|\yy'(x)|>50$ means that the tangent at $x$ makes an angle less than $1^\circ$ with the vertical line. Combined with the requirement for ${\rm Var}_J(\yy)$, the graph of $\yy$ will appear $\ort$ or $\rav$-shaped, and $x=\xi$ behaves as a vertical asymptote.


\section{The Emden-Fowler equation}\label{sct:ef}

We apply `mixed-style' techniques to investigate the {\ivp}:
\begin{align*}
u''=- x^\rr u^\pp, \;\;u(0)=1,\;u'(0)=0,\qquad(\rr>0,\;\pp>1).
\label{eq:ef}\tag{EF}
\end{align*}
Let $\uu_\exc=\uu_{\rr,\pp}$ be the exact solution. Our approach involves two stages: 
\begin{itemize}[leftmargin=3ex]
\item 
finding integrable approximations of~\eqref{eq:ef}, which yield upper/lower envelopes of the exact solution $\uu_\exc=u_{\rr,\pp}$ and analysing the relationship between them;
\item 
running numerical tests based on the analytical outcome.  
\end{itemize}
First we plot the graphs of $\uu_\exc$ and $\uu'_\exc$, for a few values of $\rr, \pp$, to get the feel of these functions and to define our objectives. 

\begin{minipage}[c]{.95\textwidth} 
	\begin{center}
		\begin{minipage}[c]{.45\textwidth}\centering 
			{\includegraphics[height=0.525\textwidth]{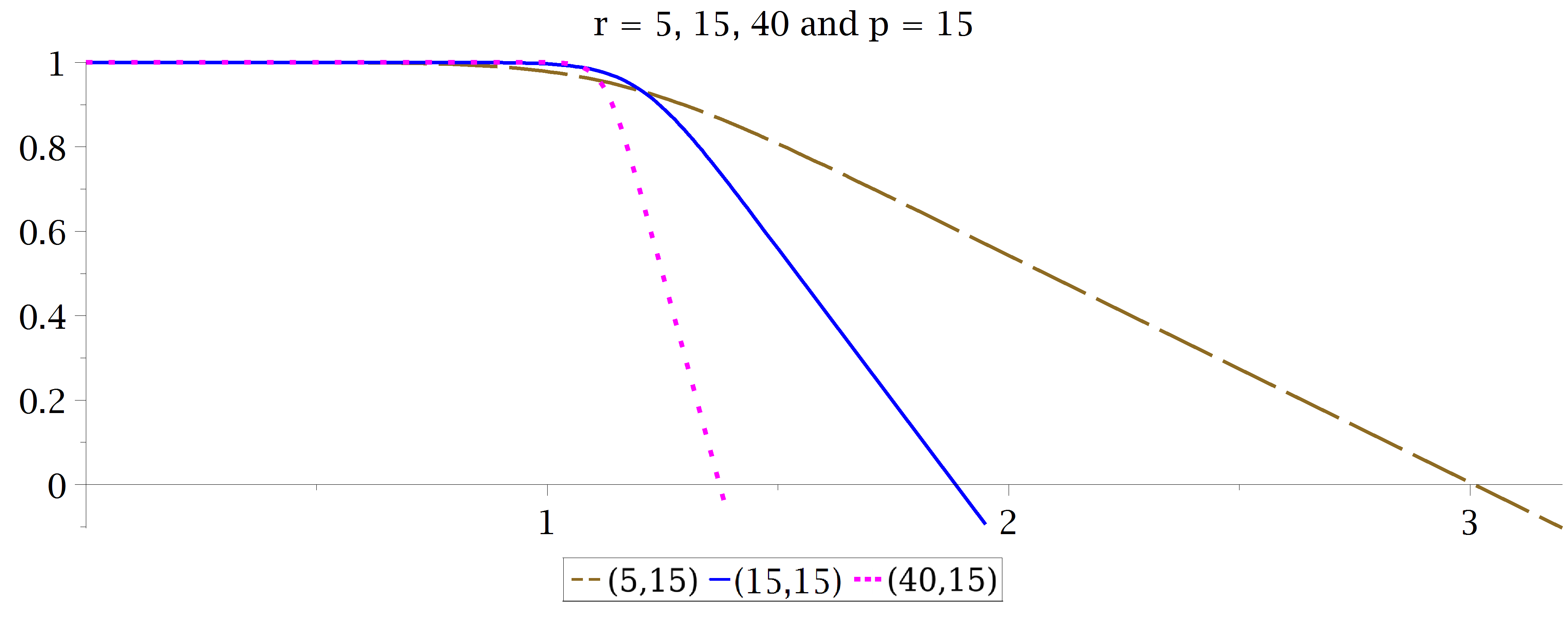}}
		\end{minipage}\hskip5ex
		\begin{minipage}[c]{.475\textwidth}\centering  
			{\includegraphics[height=0.525\textwidth]{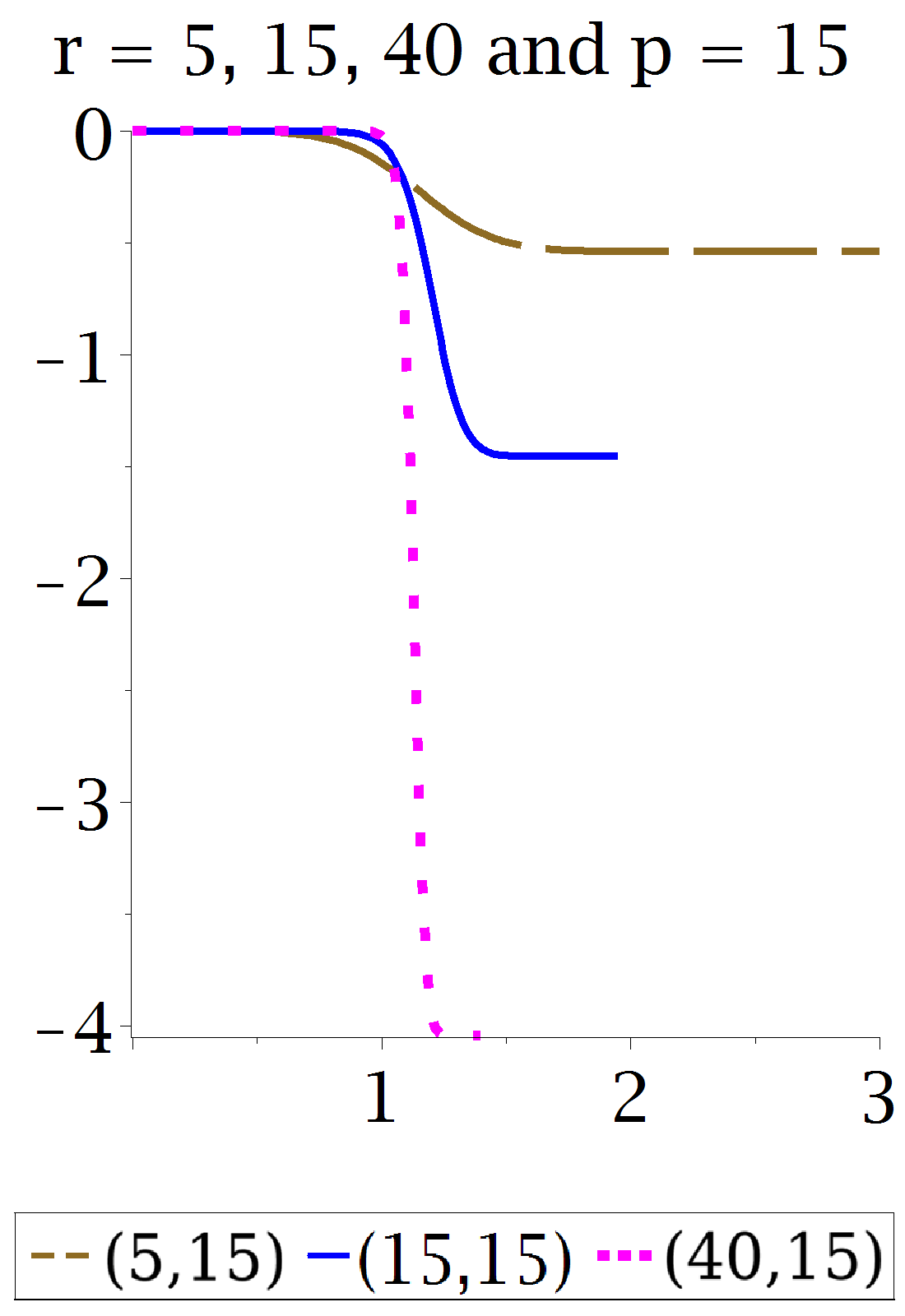}}
		\end{minipage}
		\captionof{figure}{sols.\,of\,\eqref{eq:ef} --left-- and of their derivatives --right--, varying $\rr$}
		\label{fig:vary-r}
	\end{center}
\end{minipage}

\begin{minipage}[c]{.95\textwidth} 
	\begin{center}
		\begin{minipage}[c]{.45\textwidth}\centering 
			{\includegraphics[height=0.525\textwidth]{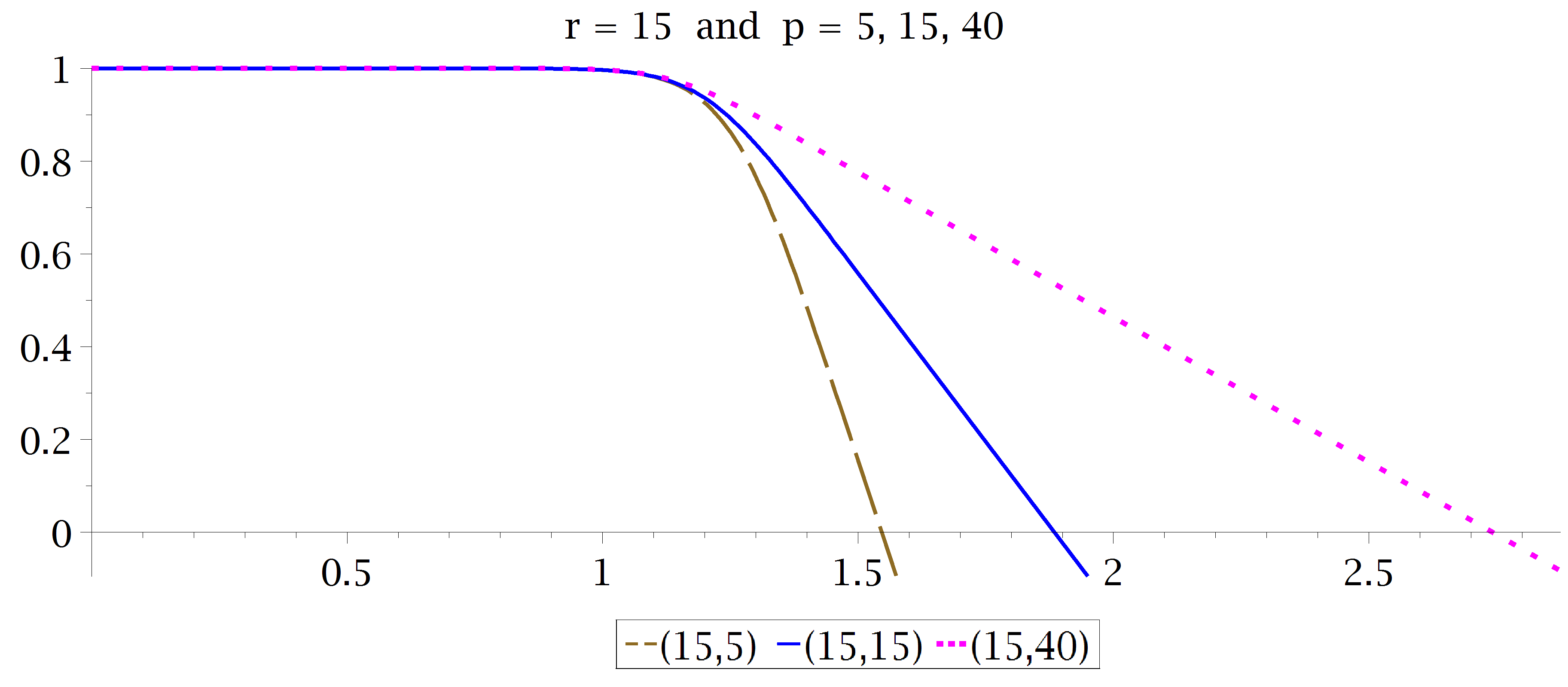}}
		\end{minipage}\hskip3ex
		\begin{minipage}[c]{.475\textwidth}\centering 
			{\includegraphics[height=0.525\textwidth]{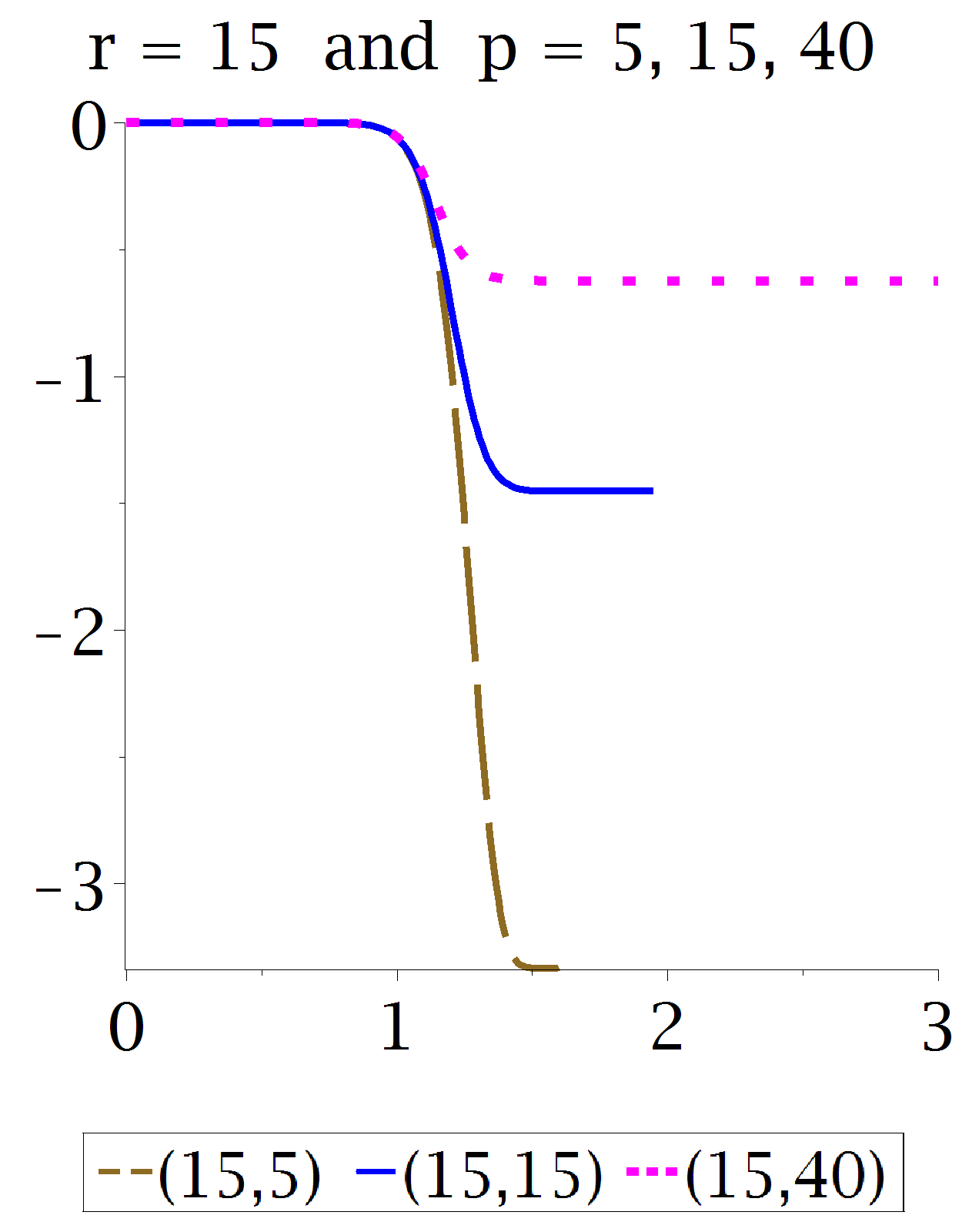}}
		\end{minipage}
		\captionof{figure}{sols.\,of\,\eqref{eq:ef} --left-- and of their derivatives --right--, varying $\pp$} \label{fig:vary-m}
	\end{center}
\end{minipage}

\nit We will consider \emph{only} the interval $[0,x_0]$, which stretches to the first zero of $u_{\rr,\pp}$. On this interval, the function has the following properties: 
\begin{enumerate}[leftmargin=5ex] 
\item  It is concave (the second derivative is negative); 
\item Strictly decreases from one until it hits zero at $x_0=x_0(\rr,\pp)$.
\item[] For irrational $\pp$, the (real) solution doesn't extend for $x>x_0$. (For $\pp\in\mbb Q$, \emph{it extends} (cf.~\cite{fow1,fow2}), and it's behaviour depends on the parity of its nominator/denominator.) 
\item Roughly, the graph of $u_\exc$ is akin to an obtuse angle. It contains two (almost) straight parts connected by a rounded region.  This is confirmed by the graphs of the derivative, similar to step functions: two (almost) horizontal steps joined by a vertical line.

\item Fig.~\ref{fig:vary-r} shows that, by increasing $\rr$ while keeping $\pp$ fixed, the graph becomes $\rav$-shaped, the tail falls steeply. On the derivative side, there is a very narrow interval between the steps, and the second step sinks quickly (the derivative changes rapidly). 
\item[] Fig.~\ref{fig:vary-m} shows that, by increasing $\pp$ while keeping $\rr$ fixed, the tail raises, and the second horizontal step raises (the derivative approaches zero). 
\end{enumerate}

\nit The equation has physical origins, so these observations raise \emph{justified, basic questions}: 
\begin{enumerate}[leftmargin=5ex]
\item[q1:] 
Can one estimate $u_{\rr,\pp}$, explicitly in terms of $\rr, \pp$ (not case-by-case)? 
\item[q2:] 
Where ends the horizontal plateau-region of the graph and starts the decreasing-part? (From a physical perspective, this might represent the end of some equilibrium state.)
\item[q3:] 
Can one estimate the (first) zero of $\uu_{\rr,\pp}$? (It represents the radius of a star.) 
\end{enumerate}
With this in mind, we employ analytical techniques to address the following objectives:
\begin{enumerate}[leftmargin=5.5ex]
\item[\S\ref{ssct:approx-fct}] 
Determine explicit envelopes which squeeze in between the exact solution $\uu_{\rr,\pp}$. Approximations for the latter are given in~\cite{ram}. The formulae in here are simpler and, in addition, we differentiate between upper/lower envelopes, which is essential for estimates. 
\item[\S\ref{ssct:turn}] 
Determine the `downward turning point' of the graph. 
\item[\S\ref{ssct:slope+bdry}] 
What condition on $\rr,\pp$ ensures the boundary layer phenomenon for $\uu_{\rr, \pp}$. 
\item[\S\ref{ssct:x0}] 
Determine the value $x_0(\rr,\pp)$ where $u_{\rr,\pp}$ vanishes. Numerically, this was investigated in~\cite{neural,wwg}, but the methods give no insight into the dependence on parameters. 
\item[\S\ref{ssct:backw}] Solve~\eqref{eq:ef} backwards, starting from $x=x_0$, to recover the initial values at $x=0$.
\item[\S\ref{ssct:zx0}] 
Given $z>1$, what relationship satisfy $\rr,\pp$, so that~\eqref{eq:ef}, with $\uu(z)=0$, is solvable.
\end{enumerate}
The $2^{\rm nd}, 3^{\rm rd}, 5^{\rm th}, 6^{\rm th}$ items above seem to have been overlooked (to our knowledge), yet they play an essential role in understanding the behaviour of the solution.


\subsection{Approximations}\label{ssct:approx-fct} 

\subsubsection{Integrable ODE}\label{sssct:ode} 

The expression $x^{\rr}\uu^\pp$ reminds the derivative of a composed function. With this motivation, we consider: 
$\;\yy''=\kp^2\yy^\nn,\;\yy(0)=1\;(\kp>0).$
Its general (decreasing) solution satisfies $\yy'=-\sqrt{2}\kp\sqrt{\frac{\yy^{\nn+1}}{\nn+1}+\CC}$, where $\CC$ is determined by $\yy'(0)$. We let $\CC=0$ and obtain $\yy(x)=\Big[1+\frac{\kp(\nn-1)}{\sqrt{2(\nn+1)}}x\Big]^{-\frac{2}{\nn-1}}$. The function 
\\[.5ex]\centerline{$
\begin{array}{l}
\YY:[0,+\infty)\to(0,1],\quad 
\YY(x):=\yy(x^{\rr+2})={\Bigl[ 1+\frac{\kp(\nn-1)}{\sqrt{2(\nn+1)}}x^{\rr+2} \Bigr]}^{-\frac{2}{\nn-1}} 
\;\text{satisfies}
\\[2ex] 
\YY''=-\kp\frac{(\rr+2)(\rr+1)\sqrt{2}}{\sqrt{\nn+1}}\,x^\rr\YY^{\frac{\nn+1}{2}}\!\cdot\!\Bigl[1-\kp\frac{(\rr+2)\sqrt{\nn+1}}{(\rr+1)\sqrt{2}}\,x^{\rr+2}\YY^{\frac{\nn-1}{2}}\Bigr],
\;\;\YY(0)=1,\,\YY'(0)=0. 
\end{array}
$}


\subsubsection{Upper envelopes}\label{sssct:up} 

The substitution $\nn=2\pp-1,\; \kp=\frac{\sqrt{\nn+1}}{(\rr+2)(\rr+1)\sqrt{2}}$ makes the resemblance with~ \eqref{eq:ef} apparent. We are led to the function: 
\begin{m-eqn}{
\begin{array}{ll}
\UP(x)={\Bigl[ 1+\frac{\pp-1}{(\rr+2)(\rr+1)}x^{\rr+2} \Bigr]}^{-\frac{1}{\pp-1}}, 
& \UP''=- \Bigl[
\underbrace{1-\frac{\pp}{(\rr+1)^2}\,x^{\rr+2}\UP^{\pp-1}}_{=:\Psi(x)}
\Bigr]\cdot x^\rr\UP^\pp, 
\\[-2ex]  
\UP(0)=1,\;\UP'(0)=0.&
\end{array}
}\label{eq:UP}
\end{m-eqn}
We compute 
$\;\Psi=1-\CC_\pp(1-\UU_+^{\pp-1})=\Bigl[ \underbrace{\CC_\pp-(\CC_\pp-1)\UU_+^{1-\pp}}_{=:\Phi(x)}\Bigr] \cdot \UU_+^{\pp-1},\;\CC_\pp:=\frac{\pp(\rr+2)}{(\pp-1)(\rr+1)},$ 
which yields $\UU_+''=-\Phi(x)\cdot x^\rr U_+^{2\pp-1}.$ In other words, $\UU_+$ also satisfies an approximate \eqref{eq:ef}-equation with index $2\pp-1$ (not only $\pp$). Consequently, we let 
\begin{m-eqn}{
\begin{minipage}{.9\textwidth}\vskip-2ex
\begin{longtable*}[l]{rll}
&$\disp \mm:=\frac{\pp+1}{2}>1,$
&$\disp \up:={\Bigl[ 1+\frac{\mm-1}{(\rr+2)(\rr+1)}x^{\rr+2} \Bigr]}^{-\frac{1}{\mm-1}},$
\\[2ex]
\text{which satisfies}
&$\up''=-\Phi(x)\cdot x^\rr\up^\pp$
&$\up(0)=1,\;\up'(0)=0,$
\\[1ex]
\text{where}
&\multicolumn{2}{l}{$\Phi(x)=\CC_\mm-(\CC_\mm-1)\up^{1-\mm}
=1-\frac{\mm+\rr+1}{(\rr+1)^2(\rr+2)}x^{\rr+2}.\; \Big(\CC_\mm=\frac{\mm(\rr+2)}{(\mm-1)(\rr+1)}\Big).
$}
\end{longtable*}
\end{minipage}
}\label{eq:up}
\end{m-eqn}
The identities~\eqref{eq:UP} and~\eqref{eq:up} are our approximate \eqref{eq:ef}-equations. 

\begin{m-lemma}\label{lm:Uu}
\begin{enumerate}[leftmargin=5ex]
\item[\rm(i)] 
$\up$ (resp. $\UP$) is decreasing and changes from concave to convex at the inflection point 
$\;\xi_\infl\!=\!\Bigl[\frac{(\rr+1)^2(\rr+2)}{\mm+\rr+1}\Bigr]^{\frac{1}{r+2}}$ 
(resp. $\;\Xi_\infl\!=\!\Bigl[\frac{(\rr+1)^2(\rr+2)}{\pp+\rr+1}\Bigr]^{\frac{1}{r+2}}$).
\item[\rm(ii)] 
$\up, \UU_+$ are \emph{upper envelopes} of $\uu_\exc=\uu_{\rr,\pp}$, satisfying 
$\uu_\exc<\up<\UU_+.$ 
\item[\rm(iii)] 
Let $\aaa(x):=\Phi(x)^{\frac{1}{r+2}}.$ Then, for any $\xi\in(0,\xi_\infl)$, one has the inequalities: 
$$
\up(x)\les\uu_\exc(\aaa(\xi)\cdot x),\;\forall x\in[0,\xi],
\text{\;\;so\;\;}\up(\aaa(\xi)^{-1}x)\les\uu_\exc(x)\les\up(x),\;\;\forall\,x\in[0,\aaa(\xi)\xi].
$$  
\item[\rm(iv)] 
For $x\!<\!{\Big[\frac{(\rr+1)^2(\rr+2)}{2(\mm+\rr+1)}\Big]}^{\frac{1}{\rr+2}}\!=\!2^{-\frac{1}{\rr+2}}\xi_\infl$, let  $\aph(x)\!:=\!\frac{1-(x/\xi_\infl)^{\rr+2}} {1-2(x/\xi_\infl)^{\rr+2}}.$ 
One has the inequality: 
\\[.5ex]\centerline{
$\aph(x)\cdot\aaa(x)\cdot\uu_\exc'(\aaa(x)x)\les\up'(x).$
}
\end{enumerate}
Statements analogous to {\rm(iii)},\,{\rm(iv)} hold true for $\UP, \bb:=\Psi^{\frac{1}{\rr+2}}$ instead of $\up,\aaa$, respectively.
\end{m-lemma}
The last two statements should be viewed as follows. About any $\xi$ slightly smaller than the inflection point, the function $\uu_\exc(x)$ behaves as $\up(\gma x)$, for some `\emph{correction factor}' $\gma=O(\aaa(\xi)^{-1})$. The estimates are useful when $\aaa(\xi)$ is close to one (note $\aaa(\xi)<1$). 

\begin{m-proof} 
(i) The statement follows from~\eqref{eq:up}. 

\nit(ii) We apply Proposition~\ref{prop:+-}(ii). Take $\veps>0$ small, let $\uu_{(\veps)}$ be defined the same as $\up$, with $(1-\veps)+\dots$ instead of $1+\dots$. It satisfies an ODE similar~\eqref{eq:up}, with $\Phi$-factor less than one, $\uu_{(\veps)}(0)>1, \uu_{(\veps)}'(0)=0$. Thus $\uu_{(\veps)}\ges\uu_\exc$, for all $x$, so we can take the limit $\veps\to0^+$. 

\nit(iii) Let $a=\aaa(\xi)$. The function $\uu_a(x):=\uu_\exc(ax)$ satisfies the ODE $\uu_a''=-a^{\rr+2}x^\rr \uu_a^\pp$, while $\up''=-\Phi(x)\cdot x^r\up^\pp\les -a^{\rr+2}x^\rr \up^\pp$. We conclude as before. 

\nit(iv) The following holds true:
\\[.5ex]\centerline{
$\begin{array}{rl}
\up'(x)
&
=\int_0^{x}\up''(t)\rd t=-\int_0^x\aaa(t)^{\rr+2}t^\rr\up(t)^\pp\rd t 
\srel{\rm(iii)}{\ges}-\int_0^x\aaa(t)^{\rr+2}t^\rr\uu_\exc\big(\aaa(t)t\big)^\pp\rd t
\\[1.5ex] 
&
=\int_0^x \aaa(t)^2\uu_\exc''\big(\aaa(t)t\big)\rd t
=\int_0^x\frac{\aaa(t)^2}{\aaa(t)+t\aaa'(t)}\cdot \frac{\rd}{\rd t}\big(\uu_\exc'\big(\aaa(t)t\big)\big)\rd t.
\end{array}$
}\\[.5ex]
A computation shows that $\frac{\aaa(t)^2}{\aaa(t)+t\aaa'(t)}$ is increasing, so it's smaller than 
\\[.5ex]\centerline{$\frac{\aaa(x)^2}{\aaa(x)+x\aaa'(x)}
=\aaa(x)\cdot\frac{1-(x/\xi_\infl)^{\rr+2}} {1-2(x/\xi_\infl)^{\rr+2}}.$}\\[.5ex] 
The derivative under the last integral is negative, so the conclusion follows. 
\end{m-proof}

\begin{m-remark} 
(i) The importance of the inflection point relies in the fact that the largest value $x$ for which one may expect $\up(x)$ to (acceptably) approximate $u_\exc(x)$ is $\xi_\infl$. For $x>\xi_\infl$, $u_\exc$ is concave and $\up$ is convex, their graphs are heading in divergent directions. 

\nit(ii) The function $\UU_+$ reminds the first term of  Ramnath's~\cite[eq. (3.2.24)]{ram} approximate solution ${\Bigl[ b_0+\frac{(\pp-1)x^{\rr+2}}{(\rr+1)(\rr+2)}-b_1x \Bigr]}^{\frac{-1}{\pp-1}}$. Our simplified formula is essential for explicit computations. 

\nit(iii) We have two upper envelopes $\up<\UP$. Although the former is sharper, we will also use the latter because it leads to formulae (involving $\rr, \pp$) which are easier to interpret. As a general rule, $\up$ will be used to approximate $\uu_\exc$ itself, especially for numerical purposes, but some key quantities (turn-points, slopes) will be estimated using $\UP$. 
\end{m-remark}

\subsubsection{Lower envelope}\label{sssct:down} 

Note that $\uu_\exc=\uu_{\rr, \pp}$ is (usually) defined on $[0,x_0]$, while $\UP, \up$ are defined on $\mbb R$ and don't vanish, so there is no similar global lower envelope. To construct one, note that upper envelopes of $u_\exc$ yield, by integration, lower envelopes, due to the `$-$' sign in~\eqref{eq:ef}. 
Since $\uu(x)\les1,$ we start with the simplest choice $x\srel{\bone}{\lmt} 1$: 
\begin{m-eqn}{
\uu_\exc(x)\ges\uu_\blt(x):=1-\frac{x^{r+2}}{(r+1)(r+2)}. 
}\label{eq:u-blt}
\end{m-eqn}
The right-hand side is positive only for $0\les x\les \xi_\blt:=[(r+1)(r+2)]^{1/(r+2)}$; the value of $\xi_\blt$ is typically close to $1$, and it is greater than $\xi_\infl$.  
\begin{itemize}[leftmargin=3ex] 
\item 
This explains the plateau regions in Fig.~\ref{fig:vary-r},~\ref{fig:vary-m}, stretching slightly further than $x=1$: it's because $x^{\rr+2}\approx0$ for $x$ in most of $[0,1]$, especially for larger values of $\rr$.
\item 
Since $\uu_\exc(\xi_\blt)>0$, one has $x_0>\xi_\blt$, so $\uu_\exc$ is always defined on $[0,\xi_\blt]$.
\end{itemize}

To improve $\uu_\blt$, the next choice would be either $\up$ or $\UU_+$, instead of $\bone$, but this leads to non-integrable expressions. The only option is finding an easily integrable function between $\UP^\pp$ and $\bone$, which approximates $\UP^\pp$ well, too. This leads to the truncated power series expansion of $\UP^\pp$ (\emph{away} from the origin, as $x_0$ can be large; near-$0$ expansions are used in \S\ref{ssct:smallp}): 
\\[.75ex]\centerline{$
f:[0,+\infty)\to\mbb R,\;\;f(t)=
\left\{
\begin{array}{cl}
1,&0\les t\les T_{\rr,\pp}:={\bigl[\frac{(\rr+1)(\rr+2)}{\pp-1}\bigr]}^{\frac{1}{r+2}},
\\[0ex] 
{\Bigl[ \frac{(\rr+1)(\rr+2)}{(\pp-1)t^{\rr+2}} \Bigr]}^{\frac{\pp}{\pp-1}},&\qquad t>T_{\rr,\pp}.
\end{array}
\right.
$}\\[.75ex]  
(The reason for choosing $\UP$ instead of $\up$ is that the truncation by $1$ is shorter.) 
By integrating $\UM''=-t^\rr f(t)$, we obtain the $\euf C^2$-function $\UM:[0,+\infty)\to\mbb R$, which provides a lower envelope for $\uu_\exc$ and for its derivative. 
\\\kern-1ex\begin{tabular}{rl}
\multicolumn{2}{l}{$-$\;\;The expression of $\UM'(x)$ is:}
\\ 
\scalebox{.85}{$\disp0\les x\les T_{\rr,\pp}:\;\UM'(x)=-\frac{x^{\rr+1}}{\rr+1}$};
&\;
\scalebox{.85}{$\disp x>T_{\rr,\pp}:\;\;\UM'(x)=-\frac{{(T_{\rr,\pp})}^{\rr+1}(\rr+2)\pp}{(\rr+1)(\rr+\pp+1)}+\frac{{(T_{\rr,\pp})}^{\frac{\pp(\rr+2)}{\pp-1}}(\pp-1)}{\rr+\pp+1}\cdot x^{-\frac{\rr+\pp+1}{\pp-1}}$.}
\end{tabular}
\\\kern-1ex\begin{tabular}{rl}
\multicolumn{2}{l}{$-$\;\;The expression of $\UM(x)$ is:}
\\[1ex]
\scalebox{.85}{$\disp0\les x\les T_{\rr,\pp}:$}&
\scalebox{.85}{$\disp\;\UM(x)=1-\frac{x^{\rr+2}}{(\rr+1)(\rr+2)}$;}
\\[1ex] 
\scalebox{.85}{$x>T_{\rr,\pp}:$}&
\scalebox{.85}{$\disp\;\UM(x)=\frac{(\rr+2)\pp-1}{\pp-1}-\frac{{(\rr+2)}^{\frac{2\rr+3}{\rr+2}}\frac{\pp}{{(\pp-1)}^{\frac{\rr+1}{\rr+2}}}}{{(\rr+1)^\frac{1}{\rr+2}}(\rr+\pp+1)}\cdot x-\frac{{(\rr+1)}^{\frac{\pp}{\pp-1}}{(\rr+2)}^{\frac{1}{\pp-1}}{(\pp-1)}^{\frac{\pp-2}{\pp-1}}}{\rr+\pp+1}\cdot x^{-\frac{\rr+2}{\pp-1}}.$}
\end{tabular}\\[1ex]
\nit In spite of its unpleasant form, $\UM$ is explicit and has the same monotonicity/geometric properties as $\uu_\exc$: it is strictly decreasing, concave, so it reaches zero before $x_0$. 

\begin{minipage}[c]{.9\textwidth}\centering
{\includegraphics[height=0.24\textwidth]{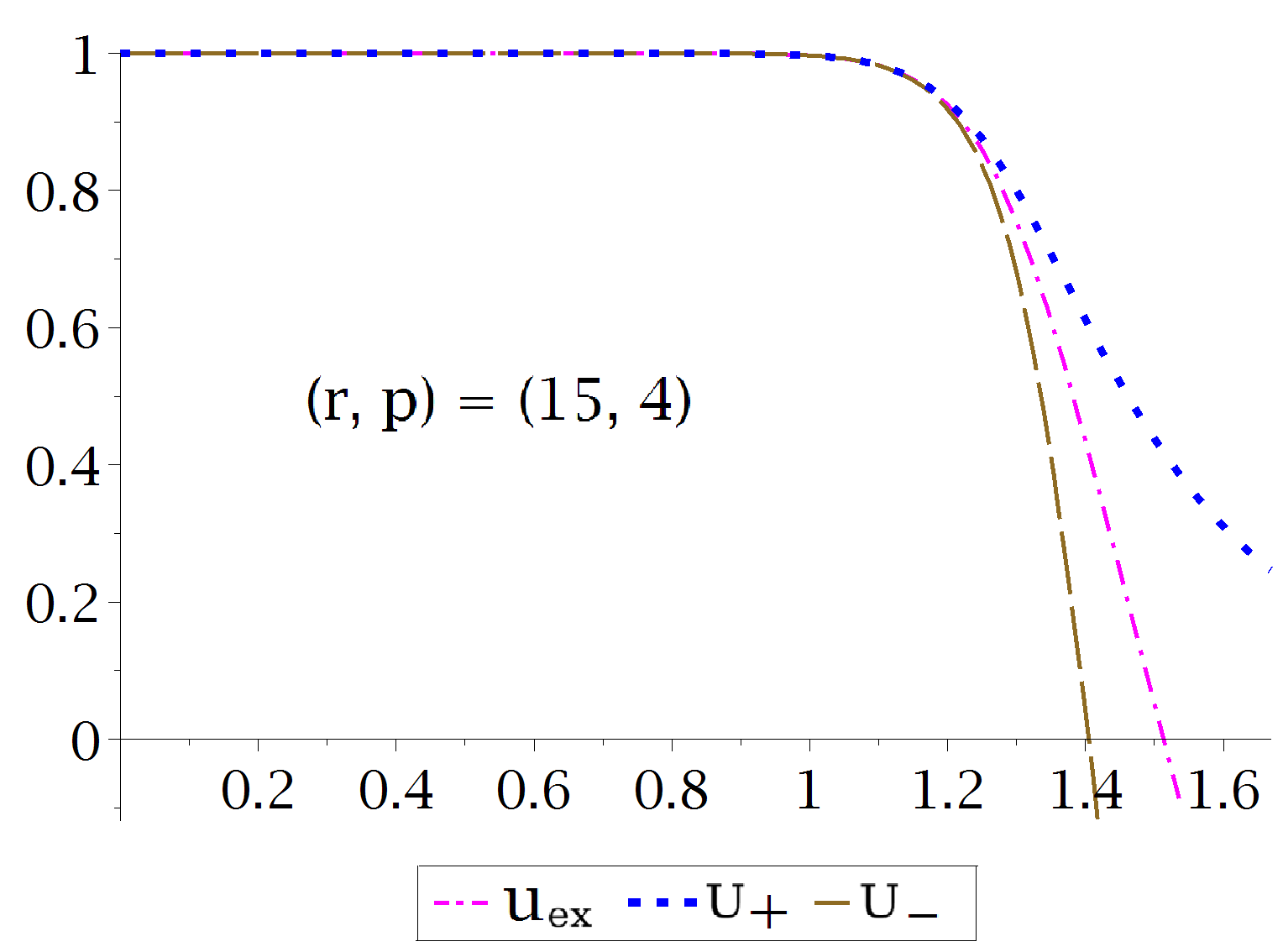}}
\hskip3ex
{\includegraphics[height=0.24\textwidth]{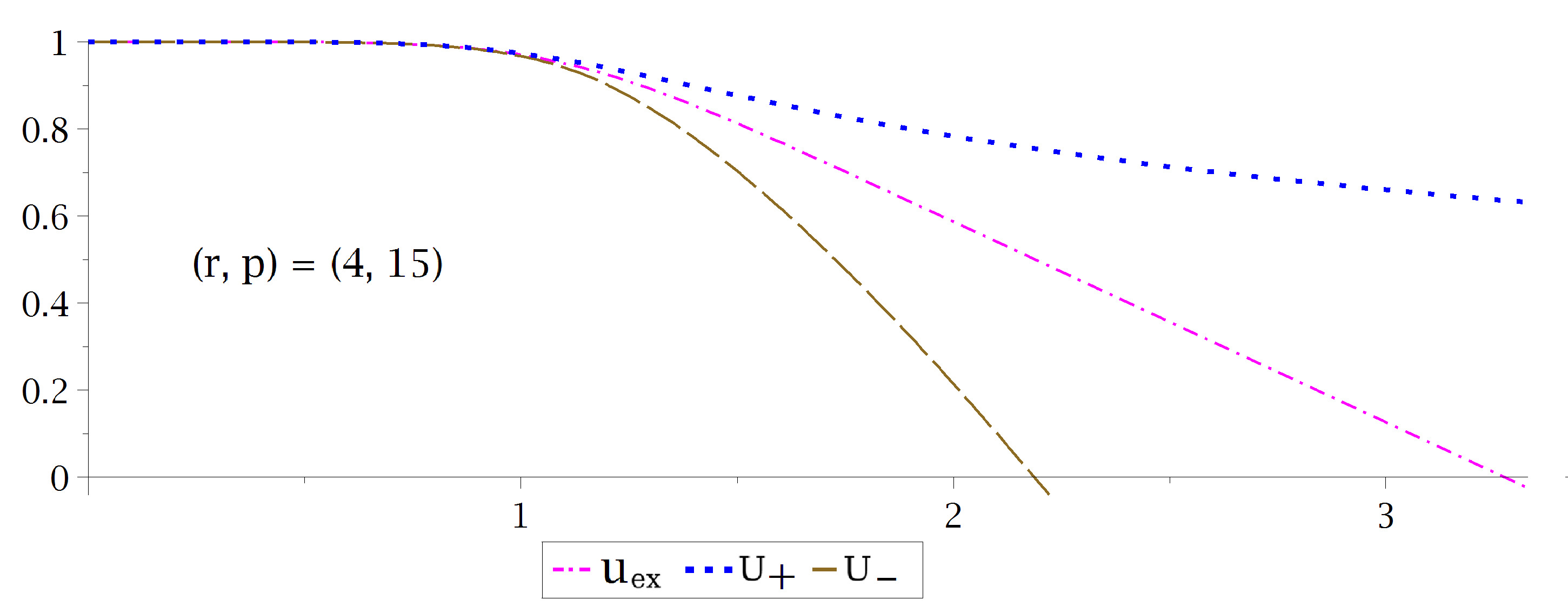}}
\captionof{figure}{solution of~\eqref{eq:ef}, $\UM$, and $\UP$}\label{fig:odeupum}
\end{minipage}

\nit The figure reveals that the graphs seem to coincide ---the envelopes are good--- only up to the `turning point', to be defined, indicating its importance. (See numerical data below.) 
\renewcommand{\arraystretch}{1.1}\textscale{.75}{ 
\begin{longtable}[c]{rl|c||ll||ll||ll||ll|}
\cline{3-11}
 &  & \diagbox[width=6ex, height=5ex]{$\pp$}{$\rr$} 
 & \multicolumn{2}{c||}{4} & \multicolumn{2}{c||}{60} 
 & \multicolumn{2}{c||}{500} & \multicolumn{2}{c|}{1700} 
\\ \hhline{--=========}
\multicolumn{1}{|r|}{$\Xi^-_\turn$} & \multicolumn{1}{|l|}{$\Delta_\pm(\Xi^-_\turn)$}  
& \multirow{2}{*}{4} & \multicolumn{1}{l|}{1.122} & 7.70(-3) 
& \multicolumn{1}{l|}{1.066} & 4.05(-4) 
& \multicolumn{1}{l|}{1.012} & 7.68(-6) 
& \multicolumn{1}{l|}{1.004} & 6.84(-7) 
\\ \cline{1-2} \cline{4-11} 
\multicolumn{1}{|r|}{$0.75\Xi^-_\turn$} & \multicolumn{1}{|l|}{$\Delta_\pm(0.75\Xi^-_\turn)$}  
&  & \multicolumn{1}{l|}{0.842} & 2.74(-4) 
& \multicolumn{1}{l|}{0.800} & 1.35(-19) 
& \multicolumn{1}{l|}{0.759} & 2.81(-131) 
& \multicolumn{1}{l|}{0.753} & 3.50(-432) 
\\ \cline{1-2}\hhline{~~=========}
 &  & \multirow{2}{*}{40} 
 & \multicolumn{1}{l|}{0.811} & 1.46(-3) 
 & \multicolumn{1}{l|}{1.054} & 8.36(-4) 
 & \multicolumn{1}{l|}{1.012} & 5.67(-5) 
 & \multicolumn{1}{l|}{1.004} & 6.21(-6) 
 \\ \cline{4-11} 
 &  &  
 & \multicolumn{1}{l|}{0.608} & 5.50(-5) 
 & \multicolumn{1}{l|}{0.791} & 3.18(-19) 
 & \multicolumn{1}{l|}{0.759} & 2.16(-130) 
 & \multicolumn{1}{l|}{0.753} & 3.22(-431) 
 \\ \hhline{~~=========} 
 &  & \multirow{2}{*}{100} 
 & \multicolumn{1}{l|}{0.699} & 6.14(-4) 
 & \multicolumn{1}{l|}{1.043} & 5.75(-4) 
 & \multicolumn{1}{l|}{1.011} & 9.29(-5) 
 & \multicolumn{1}{l|}{1.004} & 1.34(-5) 
 \\ \cline{4-11} 
 &  &  
 & \multicolumn{1}{l|}{0.524} & 2.33(-5) 
 & \multicolumn{1}{l|}{0.782} & 2.30(-19) 
 & \multicolumn{1}{l|}{0.758} & 3.70(-130) 
 & \multicolumn{1}{l|}{0.753} & 7.07(-431) 
 \\  \hhline{~~=========} 
 &  & \multirow{2}{*}{1000} 
 & \multicolumn{1}{l|}{0.478} & 6.33(-5) 
 & \multicolumn{1}{l|}{1.009} & 8.71(-5) 
 & \multicolumn{1}{l|}{1.009} & 6.33(-5) 
 & \multicolumn{1}{l|}{1.003} & 3.10(-5) 
 \\  \cline{4-11} 
 &  &  
 & \multicolumn{1}{l|}{0.359} & 2.41(-6) 
 & \multicolumn{1}{l|}{0.757} & 3.67(-20) 
 & \multicolumn{1}{l|}{0.757} & 2.90(-130) 
 & \multicolumn{1}{l|}{0.752} & 1.86(-430) 
\\ \cline{3-11} 
\caption{Error analysis for $\Delta_\pm=\UP-\UM$ about $\Xi^-_\turn$ \eqref{eq:x-turn}.}\label{tab:error}
\end{longtable}
}\renewcommand{\arraystretch}{1} 
\nit Numerical estimates of $\uu_\exc$ on longer intervals usually involve power series. This could be done here, too, by replacing $f$ with higher order expansions. But for our purposes ---extracting analytical information from approximating functions--- this method is of no use. We take a completely different approach: the descending part of the graph of $\uu_\exc$ is almost straight (see \S\ref{ssct:slope+bdry}), so Newton's method at the turning point will approximate $\uu_\exc$ by a line.


\subsection{Turning point}\label{ssct:turn}

First we clarify its defining property: $u''_{\rr, \pp}(x)$ vanishes at $x=0$ and at $x_0(\rr, \pp)$, it's strictly negative in between. Thus it has an absolute minimum at $x=\xi_\turn$ where $u_{\rr,\pp}'''(x)=\bigl( x^\rr\uu_{\rr, \pp}(x)^\pp \bigr)'=0$, so  
\begin{m-eqn}{
x=\xi_\turn\;\;\text{solves}\quad 
\rr\cdot u_{\rr,\pp}(x)+\pp\cdot x u_{\rr,\pp}'(x)=0.
}\label{eq:turn}
\end{m-eqn}
This is indeed a turning point in the na\"ive sense: the graph is the most bended downward. The equation can't be solved exactly, so we approximate it:  
\begin{longtable*}[l]{rcl}
$(x^\rr\cdot\UP(x)^\pp)'=0$
&$\;\Rightarrow\;$& 
$x=\Xi^+_\turn:=\Bigl[\frac{\rr(\rr+1)(\rr+2)}{\rr+2\pp}\Bigr]^{\frac{1}{\rr+2}};$
\\[.5ex]
$(x^\rr\cdot\up(x)^\pp)'=0$
&$\;\Rightarrow\;$& 
$x=\xi^+_\turn:=\Bigl[\frac{\rr(\rr+1)(\rr+2)}{\rr\mm+2(2\mm-1)}\Bigr]^{\frac{1}{\rr+2}}
=\Bigl[\frac{\rr(\rr+1)(\rr+2)}{\rr\frac{\pp+1}{2}+2\pp}\Bigr]^{\frac{1}{\rr+2}}<\Xi^+_\turn.$
\\[.5ex]
\multicolumn{3}{l}{Note that $\xi^+_\turn<\Bigl[ \frac{\rr(\rr+1)(\rr+2)}{2(\rr+\mm+1)}\Bigr]^{\frac{1}{\rr+2}}<\Xi^+_\turn$.
\;(The middle term appears in Lemma~\ref{lm:Uu}(iv).)}
\end{longtable*}
\nit Since $\Psi(\Xi^+_\turn)=1/(\rr+1)$, it follows that 
$
\UP(x)\les\uu_\exc\Big( [\frac{1}{\rr+1}]^{\frac{1}{\rr+2}}\cdot x \Big)
,$ for $x\in[0,\Xi^+_\turn]$. 
But $\xi_\turn,\Xi^+_\turn$ are local maxima for $x^\rr\uu_\exc^\pp, x^\rr\UP^\pp$, respectively, so 
\begin{m-eqn}{
\Xi^-_\turn:=\Bigr[\frac{\rr(\rr+2)}{\rr+2\pp}\Bigr]^{\frac{1}{r+2}}
=\frac{\Xi^+_\turn}{(\rr+1)^{1/(\rr+2)}}
<\xi_\turn<
\Xi^+_\turn=\Bigl[\frac{\rr(\rr+1)(\rr+2)}{\rr+2\pp}\Bigr]^{\frac{1}{r+2}}.
}\label{eq:x-turn}
\end{m-eqn}
Indeed, for the first inequality, note that 
\\[.5ex]\centerline{
$\rr\uu_\exc(\Xi^-_\turn)+\pp\Xi^-_\turn\uu_\exc(\Xi^-_\turn)
>\rr\UM(\Xi^-_\turn)+\pp\Xi^-_\turn\UM(\Xi^-_\turn)=\rr\big[1-\frac{2}{\rr+2\pp}\big]>0$,
}\\[.5ex] 
so $x^\rr\uu_\exc(x)^p$ is still decreasing at $x=\Xi^-_\turn$.

Numerically, one verifies~\eqref{eq:x-turn} by plotting $r\cdot u_{r,m}(x)+m\cdot x u_{r,m}'(x)$: it changes sign (positive to negative) between $\Xi^\mp_\turn$, so this gives a method to determine $\xi_\turn$. 
Table~\ref{tab:turn} lists its values for a wide range of parameters $\rr,\pp$. Note that $\xi^+_\turn$ approximates $\xi_\turn$ better, reflecting that $\up$ is closer to $\uu_\exc$ than $\UP$. The `more generous' interval $(\Xi^-_\turn,\Xi^+_\turn)$ will be useful later on, for estimating the zero of $\uu_\exc$. 
\renewcommand{\arraystretch}{1.1}{\textscale{.75}{
\begin{longtable}[c]{cc|c||c|l||c|l||c|l||c|l||c|l||c|l|}
\hhline{~~-------------}
&&
	\diagbox[width=6ex, height=5ex]{$\pp$}{$\rr$}          
	& \multicolumn{2}{|c||}{\textbf{4}}    & \multicolumn{2}{c||}{\textbf{15}}   
	& \multicolumn{2}{c||}{\textbf{60}}   & \multicolumn{2}{c||}{\textbf{225}}  
	& \multicolumn{2}{c||}{\textbf{500}}  & \multicolumn{2}{c|}{\textbf{1700}}   
\\ \hhline{--=============}
\multicolumn{1}{|c|}{$\boldsymbol\xi_\turn$} &$\xi^+_\turn$&
	\multirow{2}{*}{\textbf{4}}   
	& \multicolumn{1}{c|}{\textbf{1.336}} & 1.371
	& \multicolumn{1}{c|}{\textbf{1.287}} & 1.302 
	& \multicolumn{1}{c|}{\textbf{1.120}} & 1.124
	& \multicolumn{1}{c|}{\textbf{1.043}} & 1.044 
	& \multicolumn{1}{c|}{\textbf{1.022}} & 1.023 
	& \multicolumn{1}{c|}{\textbf{1.0080}} & 1.0082         
\\ \cline{1-2}\cline{4-15} 
\multicolumn{1}{|c|}{$\Xi^-_\turn$}&\multicolumn{1}{|c|}{$\Xi^+_\turn$}&
    & {1.122} & 1.467 & {1.152} & 1.356 & {1.066} & 1.139 
	& {1.024} & 1.048 & {1.012} & 1.025 & {1.0043} & 1.0087 
\\ \cline{1-2}\hhline{~~=============}
&&
	\multirow{2}{*}{\textbf{15}}   
	& \multicolumn{1}{c|}{\textbf{1.101}} & 1.116
	& \multicolumn{1}{c|}{\textbf{1.208}} & 1.214
	& \multicolumn{1}{c|}{\textbf{1.101}} & 1.103
	& \multicolumn{1}{c|}{\textbf{1.038}} & 1.039
	& \multicolumn{1}{c|}{\textbf{1.020}} & 1.020
	& \multicolumn{1}{c|}{\textbf{1.0075}} & 1.0075
\\ \cline{4-15} 
&&
    & {0.943} & 1.233 & {1.107} & 1.303 & {1.061} & 1.134 
	& {1.023} & 1.048 & {1.012} & 1.024 & {1.0043} & 1.0087 
\\ \hhline{~~=============}
&&
	\multirow{2}{*}{\textbf{40}}  
	& \multicolumn{1}{c|}{\textbf{0.941}} & 0.951
	& \multicolumn{1}{c|}{\textbf{1.145}} & 1.148
	& \multicolumn{1}{c|}{\textbf{1.086}} & 1.086
	& \multicolumn{1}{c|}{\textbf{1.034}} & 1.034
	& \multicolumn{1}{c|}{\textbf{1.018}} & 1.018
	& \multicolumn{1}{c|}{\textbf{1.0069}} & 1.0069
\\ \cline{4-15} 
&&
    & {0.811} & 1.061 & {1.059} & 1.247 & {1.050} & 1.120 
	& {1.022} & 1.047 & {1.012} & 1.024 & {1.0043} & 1.0087 
\\ \hhline{~~=============} 
&&
	\multirow{2}{*}{\textbf{100}} 
	& \multicolumn{1}{c|}{\textbf{0.809}} & 0.817
	& \multicolumn{1}{c|}{\textbf{1.086}} & 1.089
	& \multicolumn{1}{c|}{\textbf{1.070}} & 1.071
	& \multicolumn{1}{c|}{\textbf{1.030}} & 1.031
	& \multicolumn{1}{c|}{\textbf{1.017}} & 1.017
	& \multicolumn{1}{c|}{\textbf{1.0064}} & 1.0064         
\\ \cline{4-15} 
&&
    & {0.699} & 0.915 & {1.010} & 1.189 & {1.043} & 1.115 
	& {1.021} & 1.046 & {1.011} & 1.024 & {1.0043} & 1.0087 
\\ \hhline{~~-------------}   
\caption{Turn-point and its lower/upper bounds.}
\label{tab:turn}
\end{longtable}
}}\renewcommand{\arraystretch}{1} 
\nit We observe that $\Xi^-_\turn\ges1$, for $\rr\ges\pp$, and $\Xi^+_\turn\les1$, for $\pp\ges(\rr+1)^3/2$. The cut-off point $T_{\rr,\pp}$ used to define $\UM$  satisfies: 
$\Xi^-_\turn\les T_{\rr,\pp}$; also, $T_{\rr,\pp}\les\Xi^+_\turn$ if $\rr\ges2+4/(\pp-2)$; e.g. $\rr,\pp\ges 4$. 
Overall, we identified the following $X$-values which help understanding $\uu_\exc$: 
\\[.5ex]\centerline{
$\;\Xi^-_\turn<\;\xi^+_\turn<\;\Xi^+_\turn<\;\Xi_\infl<\;\xi_\infl<\;\xi_\blt.$
}\\[.5ex]
The `bending' of $\uu_\exc$ happens near $\xi^+_\turn$. To estimate the amount of this bending compared to $\UP, \up$ (ratio of their slopes), we compute the values of the $\Psi$ and  $\Phi$ in~\eqref{eq:UP}. The `correction factors' discussed in Lemma~\ref{lm:Uu} are their $(\rr+2)^{\rm nd}$-roots, respectively. 
\begin{m-eqn}{
\begin{minipage}{.875\textwidth}
\textscale{.9}{\begin{longtable*}[l]{ll}
$\Phi(\Xi^-_\turn)=1-\frac{\rr\big(\rr+\frac{\pp+1}{2}\big)}{(\rr+1)^2(\rr+2\pp)}\approx1,$
&$\Psi(\Xi^-_\turn)=1-\frac{\rr(\rr+2)\pp}{(\rr+1)[\rr^2+(3\rr+2)\pp]}>1-\frac{\rr+2}{3\rr+2},$
\\[1ex]
$\Phi(\xi^+_\turn)=1-\frac{\rr(\rr+\mm+1)}{(\rr+1)(\rr\mm+2(2\mm-1))}
\approx1-\big(\frac{1}{\rr}+\frac{1}{\mm}\big),$
& $\Psi(\xi^+_\turn)=1-\frac{(\rr+2)\rr}{(\rr+1)\big( \rr\frac{3\pp-1}{2\pp}+2 \big)}\approx\frac{1}{3},$
\\[1ex] 
$\Phi(\Xi^+_\turn)= \frac{1}{\rr+1}+\frac{\rr\cdot(3\mm-3)}{(\rr+1)(\rr+(4\mm-2))},$
&$\Psi(\Xi^+_\turn)= \frac{1}{\rr+1},$
\\[1ex] 
$\Phi(\Xi_\infl)=\frac{\mm-1}{\rr+2\mm}>\frac{1}{\rr+2}$&$\Psi(\Xi_\infl)=0.$
\end{longtable*}} 
\end{minipage}}\label{eq:slopes-turn}
\end{m-eqn}
\nit The table shows that $\xi^+_\turn$ is too close to $\xi_\turn$, both $\up,\UP$ approximate well $\uu_\exc$. The difference becomes noticeable at $\Xi^+_\turn$, where $\UP, \up$ approximately solve  $\uu''=-\frac{1}{\rr+1}x^\rr\uu^\pp$ rather than~\eqref{eq:ef} itself. Thus we obtain the clear-cut estimate for the `correction factor' 
\\[0.5ex]\centerline{
$\gma\approx (\rr+1)^{\frac{1}{\rr+2}}.$
}\\[0.5ex] 
These computations clarify why in Fig.~\ref{fig:odeupum} the envelopes overlap $\uu_\exc$ on $[0, \Xi^-_\turn]$. Since $\Psi(\Xi^-_\turn)\ges 0.6$, the correction factor between $\UP$ and $\uu_\exc$ is at most $(5/3)^{1/(\rr+2)}$; for $\rr=3$, it's already less than $1.1$. Therefore $\UP$ approximates $\uu_\exc$ well on $[0,\Xi^-_\turn]$. In the boundary layer situation ($\rr/\pp>50$), $\Psi(\Xi^-_\turn)>0.98$, so the correction factor is less than $1.005$.


\subsection{Slope and boundary layer}\label{ssct:slope+bdry}

Estimating the slope of $\uu_\exc$ about $\xi_\turn$ is required by Newton's method, in order to approximate $x_0=x_0(\rr,\pp)$. For this, it's useful analysing the graph of $\uu_\exc'$ on $[\xi_\turn,x_0]$. The third derivative is $u_\exc'''=-x^{\rr-1}\uu_\exc^{\pp-1}(\rr\uu_\exc+\pp x\uu_\exc')$, so $\uu_\exc'$ is convex for $x>\xi_\turn$, its graph lies above the tangent lines. 
\\\begin{tabular}{l}
\begin{minipage}[l]{.65\textwidth}
The tangent  to $\uu_\exc'$ at $x_0$ is horizontal. 
\\\centerline{$(\uu_\exc')'(x_0)=\uu_\exc''(x_0)=0,\;\text{also}\;(\uu_\exc')''(x_0)=u_\exc'''(x_0)=0.$}\\    
So we have $\uu_\exc'(x)-\uu_\exc'(x_0)=O((x-x_0)^3)$. Note that~\eqref{eq:ef} actually implies $\uu_\exc'(x)-\uu_\exc'(x_0)=O\bigl((x-x_0)^{1+\lfloor{\pp}\rfloor}\bigr).$ 

This explains why the graph of $\uu_\exc$ is visually straight  (see Fig.~\ref{fig:vary-r}); for $x\to x_0$, $\uu_\exc'(x)$ converges quickly to $\uu_\exc'(x_0)$. 
\end{minipage}
\begin{minipage}[l]{.33\textwidth}
\begin{m-eqn}{
\scalebox{.85}{
\begin{tikzpicture}[baseline=(baseline point)]
\coordinate (baseline point) at (0,-1.5);
\draw[-] (-2.5,0) -- (0.25,0); 
\draw[thick,densely dotted] (0,-2.75) -- (0,.25); 
\node at (0.33,-0.25) {\mbox{\scriptsize$x_0$}}; \node at (0,0) {$\bullet$};
\draw[thick, densely dotted] (-2,-2.5) -- (-2,0.1) node[above] {\mbox{\kern2ex\scriptsize$\xi^+_\turn$}};
\node at (-2,0) {$\bullet$}; 
\node at (-2.5,-0.25) {\mbox{\scriptsize$\xi_\turn$}};\node at (-2.3,0) {$\bullet$}; 
\draw[thick,dotted] (-2.3,-1.25) -- (-2.3,0.1);
\node at (-1.33,-0.25) {\mbox{\scriptsize$\Xi_\turn$}};\node at (-1.5,0) {$\bullet$}; 
\node at (-0.75,0) {$\bullet$}; \node at (-0.75,0.35) {\mbox{\scriptsize$\Xi_\infl$}};
\draw[very thick,domain=-2.4:0] plot (\x,{1/16*\x*\x*\x*\x-2.5});
\node at (0.25,-2.25) {\mbox{\scriptsize$\boldsymbol{y=\uu'(x)}$}};
\draw[densely dashdotdotted,domain=-2.5:-1.35] plot(\x,{-1.5-2*(\x+2)});
\draw[densely dashdotdotted, thick] (-2,-2.55) -- (0,-2.55);
\end{tikzpicture}}
}\label{eq:u-der}
\end{m-eqn}
\end{minipage}
\end{tabular}


\subsubsection*{Upper bound} 

Lemma\ref{lm:Uu}(iv) implies  $\aph(\xi^+_\turn)\cdot\aaa(\xi^+_\turn)\cdot\uu'_\exc\big(\aaa(\xi^+_\turn)\xi^+_\turn\big)\les\up'(\xi^+_\turn)$; we used that $(\xi^+_\turn/\xi_\infl)^{\rr+2}\les1/2$. We compute  
$\aph(\xi^+_\turn)\approx\frac{\mm-1}{\mm-2}$ and $\aaa(\xi^+_\turn)\approx\big( \frac{\rr}{\rr+1} \big)^{\frac{1}{\rr+2}}$, therefore we declare that $\up'(\xi^+_\turn)$ is an upper bound for $\uu'_\exc(\xi^+_\turn)$; in all numerical tests, it is so.
\begin{m-eqn}{
\scalebox{1}{$
\up'(\xi^+_\turn)
=-\frac{\rr\cdot\up(\xi^+_\turn)}{\pp\cdot\xi^+_\turn}
=-\frac{\rr}{\pp}\cdot\Big[ \rr\frac{\pp+1}{2}+2\pp \Big]^{\frac{1}{\rr+2}}
\cdot\Big[ \frac{\rr\frac{\pp+1}{2\pp}+2}{\rr+2} \Big]^{\frac{1}{\pp-1}}=O\Bigl(\frac{\rr}{\pp^{1-\frac{1}{\rr+2}}}\Bigr). 
$}
}\label{eq:S+}
\end{m-eqn}
The right-hand side is strongly negative for $\rr/\pp> 50$, so $\uu_\exc$ will have boundary layer. 

\subsubsection*{Lower bound} 
 
Since $\UM'(x)$ is obtained by integrating $-t^\rr\cdot f(t)<-t^\rr\cdot\uu_\exc(t)^\pp$, we have: 
\renewcommand{\arraystretch}{1.1}{\textscale{1}{
\begin{longtable*}[l]{rl}
$\uu_\exc(\Xi^+_\turn) >$
&
$\kern-.5ex\UM'(\Xi^+_\turn)=-\frac{(\rr+2)^2\cdot \pp}{(\rr+\pp+1)\cdot (\pp-1)} {\bigl[\frac{\pp-1}{(\rr+1)(\rr+2)}\bigr]}^{\frac{1}{\rr+2}}+\frac{{[(\rr+1)(\rr+2)]}^{\frac{\rr+1}{\rr+2}}}{(\rr+\pp+1)\cdot{(\pp-1)}^{\frac{1}{\pp-1}}} {\bigl[1+\frac{2\pp}{\rr}\bigr]}^{\frac{1}{\rr+2}+\frac{1}{\pp-1}}.$
\\[1.5ex] 
$(\rr>50\pp)\leadsto$
&
$\approx  -\frac{{(\rr+2)}^2}{\Xi^+_\turn(\rr+\pp+1)}\cdot \Bigl[\frac{\pp}{\pp-1}-\frac{\rr+1}{\rr+2}\bigl( 1+\frac{2}{\rr} \bigr)\Bigr]= -\frac{{\rr+2}}{\Xi^+_\turn\cdot(\pp-1)}=O\bigl(\frac{\rr}{\pp}\bigr).$
\\[1.5ex] 
$(\pp\gg \rr)\leadsto$
&
$\approx
-\frac{\pp^{\frac{1}{\rr+2}}\cdot(\rr+2)^{\frac{\rr}{\rr+2}}}{\pp+\rr+1}\cdot(\rr+2)\bigl( 1-{\rr^{-\frac{1}{\rr+2}}} \bigr)
\approx-\frac{(\rr+2)^{1-\frac{2}{\rr+2}}}{\pp^{ 1-\frac{1}{\rr+2} }}\cdot\ln(\rr).$ 
\end{longtable*}
}

\begin{m-lemma}\label{lm:bdlayer}
\nit{\rm(i)} We have $\UM'(\Xi^+_\turn)<\uu_\exc(\Xi^+_\turn)$ and $\uu_\exc(\xi^+_\turn)<\up'(\xi^+_\turn)$.

\nit{\rm(ii)} For $\rr>50\pp$, the solution of~\eqref{eq:ef} has boundary layer. 
 (For small $\pp$, one should use~\eqref{eq:S+}.) 
The layer is located in the interval $(\xi^+_\turn,x_0)$, which contains $(\xi^+_\turn,\xi_\blt)$. 
\end{m-lemma}

For the boundary layer, we apply Definition~\ref{def:bl} to $I=[0,x_0]\supset J=[\xi^+_\turn,x_0]$. Note also that, since $\uu(1)>\uu_\blt(1)=1-[(\rr+1)(\rr+2)]^{-1}$, the (negative of the) slope of the secant joining $(1,\uu(1))$ and $(x_0,0)$ is $(x_0-1)\big/\uu(1)\big.<(x_0-1)\big/\big(1-[(\rr+1)(\rr+2)]^{-1}\big)\Big..$ 
Thus $\;-(\text{slope of secant})<{1}/{50},\,\;\text{for}\,\; x_0<1.015.$ 
So we have boundary layer whenever $x_0<1.015$. 

Below, we computed slopes at $\Xi^+_\turn$, which better reflects the behaviour of $\uu_\exc'$ (see~\eqref{eq:slopes-turn}). As $\pp=2\mm-1$ increases, the ratio $\frac{\uu_\exc'(\Xi^+_\turn)}{\up'(\Xi^+_\turn)}$ approaches $1$, in agreement with the correction factor $\Phi(\Xi^+_\turn)^{\frac{1}{\rr+2}}\srel{\mm\gg0}{\approx}[\frac{1+0.75\cdot\rr}{\rr+1}]^{\frac{1}{\rr+2}}$. 
The shaded entries correspond to boundary layer cases. The corresponding values of $\rr,\pp$ respect the inequality $\rr/\pp>50$. 
\renewcommand{\arraystretch}{1.1}{\textscale{.75}{
\begin{longtable}[c]{cc|c||r|r||r|r||r|r||r|r||r|r||r|r|}
\hhline{~~-------------}
&&
	\diagbox[width=6ex, height=5ex]{$\pp$}{$\rr$}
    & \multicolumn{2}{c||}{4} & \multicolumn{2}{c||}{15} 
	& \multicolumn{2}{c||}{60} & \multicolumn{2}{c||}{225} 
	& \multicolumn{2}{c||}{500} & \multicolumn{2}{c|}{1700} 
\\ \hhline{--=============}
\multicolumn{1}{|c|}{$\uu_\exc'(\Xi^+_\turn)$}&
\multicolumn{1}{c|}{\multirow{2}{*}{$\frac{\uu'_\exc(\Xi^+_\turn)}{\up'(\Xi^+_\turn)}$}}& 
	\multirow{2}{*}{4} 
	& {-0.82}  & \multirow{2}{*}{1.19} 
	& {-3.31}  & \multirow{2}{*}{1.26} 
	& {-15.32} & \multirow{2}{*}{1.30} 
	& \cellcolor{lgray}{\textbf{-61.73}} & \multirow{2}{*}{1.31} 
	& \cellcolor{lgray}{\textbf{-139.99}} & \multirow{2}{*}{1.31} 
	& \cellcolor{lgray}{\textbf{-482.87}} & \multirow{2}{*}{1.31} 
\\ \hhline{-~~-~-~-~-~-~-~}
	\multicolumn{1}{|c|}{$\up'(\Xi^+_\turn)$}&\multicolumn{1}{c|}{}&
    & {-0.69}  & & {-2.62}  & & {-11.77} & 
	& {-46.98} & & {-106.36} & & {-366.34} & 
\\ \cline{1-2}\hhline{~~=============} 
&&
	\multirow{2}{*}{15}   
	& {-0.31} & \multirow{2}{*}{1.10} 
	& {-1.23} & \multirow{2}{*}{1.12} 
	& {-5.73}  & \multirow{2}{*}{1.14} 
	& {-23.14} & \multirow{2}{*}{1.15} 
	& \cellcolor{lgray}{\textbf{-52.51}} & \multirow{2}{*}{1.15} 
	& \cellcolor{lgray}{\textbf{-181.20}} & \multirow{2}{*}{1.15} 
\\ \hhline{~~~-~-~-~-~-~-~}
&&
    & {-0.28}  & & {-1.09}  &  	& {-5.01}  & 
	& {-20.08} & & {-45.47} & & {-156.70} &  
\\ \hhline{~~=============} 
&&
	\multirow{2}{*}{40}   
	& {-0.14}  & \multirow{2}{*}{1.07} 
	& {-0.52}  & \multirow{2}{*}{1.07} 
	& {-2.42}  & \multirow{2}{*}{1.07} 
	& {-9.77}  & \multirow{2}{*}{1.08} 
	& {-22.18}  & \multirow{2}{*}{1.08} 
	& \cellcolor{lgray}{\textbf{-76.54}}  & \multirow{2}{*}{1.08} 
\\ \hhline{~~~-~-~-~-~-~-~}
&&
    & {-0.13}  & & {-0.49}  & 	& {-2.24}  & 
	& {-9.03}  & & {-20.45} & & {-70.48} & 
\\ \hhline{~~=============} 
&&
	\multirow{2}{*}{100}  
	& {-0.068} & \multirow{2}{*}{1.06} 
	& {-0.23}  & \multirow{2}{*}{1.05} 
	& {-1.02}  & \multirow{2}{*}{1.04} 
	& {-4.13}  & \multirow{2}{*}{1.04} 
	& {-9.36}   & \multirow{2}{*}{1.04} 
	& {-32.30}  & \multirow{2}{*}{1.04} 
\\ \hhline{~~~-~-~-~-~-~-~}
&&
    & {-0.064} & & {-0.217}  & 
	& {-0.98}  & & {-3.95}  & 
	& {-8.95}  & & {-30.86}  & 
\\ \hhline{~~=============} 
&&
	\multirow{2}{*}{1000} 
	& {-0.010} & \multirow{2}{*}{1.06} 
	& {-0.026} & \multirow{2}{*}{1.04} 
	& {-0.100} & \multirow{2}{*}{1.02} 
	& {-0.43}  & \multirow{2}{*}{1.01} 
	& {-0.97}   & \multirow{2}{*}{1.01} 
	& {-3.35}   & \multirow{2}{*}{1.009} 
\\ \hhline{~~~-~-~-~-~-~-~}
&&
    & {-0.009} & & {-0.025} & 	& {-0.107} &
	& {-0.42}  & & {-0.96}   & & {-3.32}   &   
\\ \hhline{~~-------------} 
\caption{slopes at $\Xi^+_\turn$}
\label{tab:slope}
\end{longtable}
}}\renewcommand{\arraystretch}{1}


\subsection{Approximating the zero}\label{ssct:x0}

The graph of $\uu_\exc$ is almost straight after the turning point. So we approximate it by its tangent line at an appropriate point: 
\\[.5ex]\centerline{
$y=\uu_\exc(\xi)+S\cdot(x-\xi),\;S=\text{slope}.$
}\\[.5ex]  
As explained in~\eqref{eq:u-der} and~\eqref{eq:slopes-turn}, $\xi$ has to be greater than $\xi_\turn$, actually should be slightly greater than $\xi^+_\turn$. For $\xi=\xi_\turn$, we readily obtain the following estimate for $x_0$: 
\begin{m-eqn}{
\scalebox{1}{$
x_0^\blt:=\xi_\turn+\frac{\uu_\exc(\xi_\turn)}{\uu_\exc'(\xi_\turn)}
\srel{\eqref{eq:turn}}{=}
\xi_\turn\cdot\Bigl(1+\frac{\pp}{\rr}\Bigr)\approx\xi^+_\turn\cdot\Bigl(1+\frac{\pp}{\rr}\Bigr).
$}
}\label{eq:x0turn}
\end{m-eqn}
Since $\uu_\exc$ is concave, it follows that $x_0^\blt>x_0$, it's always an upper bound. This estimate is already precise in the boundary layer situation ($\rr\ges 50\pp$), due to the almost flatness of the downhill region. The next considerations improve this: we replace (`$\mt$', for short) each term in the tangent line above with approximate upper/lower values.
\begin{itemize}[leftmargin=3ex]
\item \unbar{upper-line}:\; 
We let $\xi\mt\xi_\infl$, $\uu_\exc(\xi)\mt\up(\xi_\infl)$, and $S\mt\up'(\xi_\infl)$. 
\item[] 
The reason for this choice is that the tangent line to $\up$ at $\xi_\infl$ has the smallest $X$-intercept (denoted $x_0^+$), as $\up$ is concave for $x\les\xi_\infl$. Since $\uu_\exc\les\up$, we have $x_0<x_0^+$. 

\item \unbar{lower-line}:\; The concavity-argument doesn't work, we take an empirical approach based on the almost flatness of $\uu_\exc'$, at the right of $\xi_\turn$ (cf.~\eqref{eq:u-der}). Let: 
\item[] \hfill$\xi\mt\Xi^+_\turn,\;\;\uu_\exc(\xi)\mt\UM(\Xi^+_\turn),\;\;S\mt\UM'(\Xi^+_\turn),\;$ denote $x_0^-$ the $X$-intercept.\hfill\null 
\item[] 
One has $\UM'(\Xi^+_\turn)<\uu_\exc'(\Xi^+_\turn)$ and 
$\frac{\uu_\exc'(x_0)-\uu_\exc'(\Xi^+_\turn)}{\uu_\exc'(\Xi^+_\turn)}=\frac{\uu_\exc''(\tld x)\cdot(x_0-\Xi^+_\turn)}{\uu_\exc'(\Xi^+_\turn)}$, $\tld x\in(\Xi^+_\turn,x_0)$, is expectedly small (small numerator, large denominator), so the tangent to $\UM$ is below $\uu_\exc$. Numerically, it's true in all tests. In fact the bound is loose, so we compute $(x_0^-+x_0^+)/2$.
 
\item \unbar{mid-line}:\; We construct a point closer to $x_0$. About the turning point, the correction factor is  $\gma=(\rr+1)^{\frac{1}{\rr+2}}$ (cf.~\eqref{eq:slopes-turn}), so we substitute $S\mt(\rr+1)^{\frac{1}{\rr+2}}\up'(\xi_\infl)$. 
We take the line having this slope through $(\xi^+_\turn,\up(\xi_\infl))$, instead of $(\xi_\infl,\up(\xi_\infl))$. Its $X$-intercept is $\xi_\mdl:=\xi^+_\turn+\frac{\up(\xi_\infl)}{(\rr+1)^{\frac{1}{\rr+2}}\up'(\xi_\infl)}.$ We define the mid-point as $x_0^\mdl:=0.775\xi_\mdl+0.225x_0^+$; it oscillates about $x_0$, depending on the value of $\rr, \pp$. 
\end{itemize}

\nit The formulae for the points defined above are as follows: 
\begin{m-eqn}{\hspace{-5ex}
\begin{minipage}{0.9\textwidth}\textscale{0.87}{\nit$
\begin{array}{ll} 
x_0^+=\xi_\infl+\frac{\up(\xi_\infl)}{-\up'(\xi_\infl)}
=\Big[\frac{(\rr+1)^2(\rr+2)}{\rr+\mm+1}\Big]^{\frac{1}{\rr+2}}\big(1+\frac{\mm}{\rr+1}\big);
&\kern-1ex  
x_0^-=\Xi^+_\turn+\frac{\UM(\Xi^+_\turn)}{-\UM'(\Xi^+_\turn)};
\\[1ex] 
x_0^\mdl
=\xi_\infl\Big( \Big[0.775\frac{\xi_\mdl}{\xi_\infl}+0.225\Big]+ [0.775+0.225(\rr+1)^{\frac{1}{\rr+2}}]\frac{\mm}{(\rr+1)^{\frac{\rr+3}{\rr+2}}}\Big];
&\kern-1ex 
x_0^\pm:=(x_0^-+x_0^+)/2.
\end{array}$}
\end{minipage}}\label{eq:manyx}
\end{m-eqn}

\begin{m-lemma}\label{lm:approx-zero}
The first zero $x_0(\rr,\pp)$ of the function $u_\exc$ is situated in the interval $(x_0^{-},x_0^{+})$. 

For $\pp/\rr\les 1$, $x_0$ is approximately $1+\frac{\pp}{\rr}$, its slightly greater than one. For $\pp/\rr\gg 1$, the value of $x_0$ becomes arbitrarily large. (See data in Table~\ref{tab:x0}.)
\end{m-lemma}

\begin{m-proof}
The last statement follows from the fact that $-\uu_\exc'(\Xi^+_\turn)<-\UM'(\Xi^+_\turn)=O(\rr/\pp)$, see \S\ref{ssct:slope+bdry}, so the graph of $\uu_\exc$ is almost horizontal for $\pp\gg\rr$.
\end{m-proof}

\nit Below we listed the values of $x_0^\pm, x_0^+, x_0^\mdl$, and $x_0=x_0(\rr,\pp)$. The ratio $\rr/\pp$ covers a wide range, $0.004$--$425$. We pushed the computations as far as we could on the laptop at hand. 
\renewcommand{\arraystretch}{1.1}{\textscale{.75}{
\begin{longtable}[c]{ccc|r|r||r|r||r|r||r|r||r|r||r|r}
\hhline{~~-------------} 
&&
	\multicolumn{1}{|c||}{\diagbox[width=6ex, height=5ex]{$\pp$}{$\rr$}}
    & \multicolumn{2}{c||}{\textbf{4}} & \multicolumn{2}{c||}{\textbf{15}} 
	& \multicolumn{2}{c||}{\textbf{60}} & \multicolumn{2}{c||}{\textbf{225}} 
	& \multicolumn{2}{c||}{\textbf{500}} & \multicolumn{2}{c|}{\textbf{1700}}            
\\ \hhline{--=============} 
\multicolumn{1}{|c|}{$\bsymb{x_0(\rr,\pp)}$}&\multicolumn{1}{|c|}{$x_0^{\mdl}$}&
	\multicolumn{1}{|c||}{\multirow{2}{*}{\textbf{4}}} 
	& \textbf{2.163}  & 2.107  
	& \textbf{1.513}  & 1.510 
	& \textbf{1.170}  & 1.172
	& \cellcolor{lgray}{\textbf{1.056}} & 1.056
	& \cellcolor{lgray}{\textbf{1.028}} & 1.028 
	& \cellcolor{lgray}{\textbf{1.0097}} & \multicolumn{1}{r|}{1.0098} 
\\ \cline{1-2}\cline{4-15} 
\multicolumn{1}{|c|}{$x_0^\pm$}&\multicolumn{1}{|c|}{$x_0^+$}&
\multicolumn{1}{|c||}{}
    & {2.214}  & 2.471 
	& {1.502} & 1.594 
    & {1.166}  & 1.188 
	& {1.055} & 1.060  
    & {1.027} & 1.030  
	& {1.0095} & \multicolumn{1}{r|}{1.0102}  
\\ \cline{1-2}\hhline{~~=============} 
&&
	\multicolumn{1}{|c||}{\multirow{2}{*}{\textbf{15}}} 
	& \textbf{3.272} & 3.170 
	& \textbf{1.884} & 1.846 
	& \textbf{1.259} & 1.253 
	& \textbf{1.078} & 1.077 
	& \cellcolor{lgray}{\textbf{1.038}} & 1.037 
	& \cellcolor{lgray}{\textbf{1.0126}} & \multicolumn{1}{r|}{1.0125} 
\\ \cline{4-15} 
&&
\multicolumn{1}{|c||}{}
    & {3.412}  & 3.908  
	& {1.825} & 2.036 
    & {1.236}  & 1.289  
	& {1.072}  & 1.085  
    & {1.035} & 1.041 
	& {1.0117} & \multicolumn{1}{r|}{1.0135}  
\\ \hhline{~~=============} 
&&
	\multicolumn{1}{|c||}{\multirow{2}{*}{\textbf{40}}} 
	& \textbf{5.603} & 5.543 
	& \textbf{2.741} & 2.687  
	&\textbf{1.470} & 1.460
	& \textbf{1.132} & 1.131 
	& \textbf{1.062} & 1.061 
	& \cellcolor{lgray}{\textbf{1.0200}} & \multicolumn{1}{r|}{1.0195} 
\\ \cline{4-15} 
&&
\multicolumn{1}{|c||}{}
    & {5.953} & 6.852 
	& {2.574} & 3.022 
    & {1.397} & 1.518  
	& {1.111} & 1.143 
    & {1.052} & 1.066  
	& {1.0167} & \multicolumn{1}{r|}{1.0209}  
\\ \hhline{~~=============} 
&&
	\multicolumn{1}{|c||}{\multirow{2}{*}{\textbf{100}}} 
	& \textbf{10.586} & 10.645 
	& \textbf{4.760}  & 4.697
	& \textbf{1.985} & 1.974
	& \textbf{1.266}   & 1.264
	& \textbf{1.121} & 1.121
	& \textbf{1.036}  & \multicolumn{1}{r|}{1.036}  
\\ \cline{4-15}  
&&
\multicolumn{1}{|c||}{}
	& {11.366} & 13.100  
	& {4.339}   & 5.315 
	& {1.785}   & 2.067  
	& {1.204}   & 1.282  
	& {1.092}   & 1.128  
	& {1.028}   & \multicolumn{1}{r|}{1.038}
\\ \hhline{~~=============} 
&&
	\multicolumn{1}{|c||}{\multirow{2}{*}{\textbf{1000}}} 
	& \textbf{66.158} & 67.460 
	& \textbf{32.416} & 32.316 
	& \textbf{9.648} & 9.639
	& \textbf{3.294} & 3.294 
	& \textbf{2.027} & 2.027 
	& \textbf{1.301} & \multicolumn{1}{r|}{1.301}  
\\ \cline{4-15} 
&&
\multicolumn{1}{|c||}{}
    & {71.616}  & 82.568  
	& {28.560}  & 36.593  
    & {7.540}   & 10.143 
	& {2.594}  & 3.354  
    & {1.691}  & 2.046  
	& {1.197} & \multicolumn{1}{r|}{1.305}  
\\ \hhline{~~-------------} 
\\  \hhline{~~------------~} 
&& 
 \multicolumn{3}{|c||}{$(\rr,\pp)$}&
 \multicolumn{3}{c||}{$(10^5, 10)$} & 
 \multicolumn{3}{c||}{$(10^5, 10^5)$} & 
 \multicolumn{3}{c|}{$(10,10^5)$}& 
\\ \hhline{~~------------~} 
&& 
 \multicolumn{3}{|c||}{$(x_0,x_0^{\mdl})$ }& 
 \multicolumn{3}{c||}{$(1.000272, 1.000272$)} & 
 \multicolumn{3}{c||}{$(1.50019, 1.50021)$} & 
 \multicolumn{3}{c|}{$(2909.99, 2910.00)$} &
\\ \hhline{~~------------~}  
\caption{Lower/upper bounds for the first zero $x_0(\rr,\pp)$ of $\uu_{\rr,\pp}$.}
\label{tab:x0}
\end{longtable}
}}\renewcommand{\arraystretch}{1}


\subsection{Small \textit{p} case}\label{ssct:smallp}

So far, we focused on large values of $\pp$. But applications require also understanding the behaviour of the exact solution $\vu{\pp,\exc}$ for low values of $\pp$ (cf.~\cite{bcfg,dcbg,mon+rox,nah,ryb}), often for $|\pp-1|\approx0$. 
The solution $\vu{1,\exc}$ of~\eqref{eq:ef} for $\pp=1$ (note $\mm=1$, too) determines the behaviour of $\vu{\pp,\exc}$ for nearby values of $\pp$, due to continuous dependence on parameters.  We denote by $x_0(\pp)$ the first zero of $\vu{\pp,\exc}$. Since we deal with variable $\pp=2\mm-1$, let $\vu{\mm,+}$ be the function defined by~\eqref{eq:up}. We follow the same steps as before. 

We observe that $\vu{1,+}(x):=\exp\Big[ -\frac{x^{\rr+2}}{(\rr+1)(\rr+2)} \Big]=\ouset{\;\mm\to1}{}{\lim}\vu{\mm,+}(x)$ satisfies: 
\\[.5ex]\centerline{
$\vu{1,+}''=-\big(1-\frac{x^{\rr+2}}{(\rr+1)^2}\big)\cdot x^\rr\vu{1,+},\;\vu{1,+}(0)=1,\;\vu{1,+}'(0)=0.$
}\\[.5ex]  
By Proposition~\ref{prop:+-}, it is an upper envelope: $\vu{1,\exc}\les\vu{1,+}$. Its turn- and inflection-points are:  
\\[.5ex]\centerline{
$\xi^+_\turn(1)=\big[\rr(\rr+1)\big]^{\frac{1}{\rr+2}},\;\;\xi_\infl(1)=(\rr+1)^{\frac{2}{\rr+2}}.$
}\\[.5ex]
(They are also obtained by letting $\mm=1$ in the corresponding formulae for $\vu{1,+}$.) 

\begin{m-lemma}\label{lm:x0smallp}
The first zero of $\vu{\pp,\exc}$ satisfies ($\mm=(\pp+1)/2$): 
\\\centerline{
$\;\big[(\rr+1)(\rr+2)\big]^{\frac{1}{\rr+2}}=
\xi_\blt<x_0(\pp)\ouset{\eqref{eq:manyx}}{}{<}x_0^+(\pp)
=\Big[ \frac{(\rr+1)^2(\rr+2)}{\rr+\mm+1} \Big]^{\frac{2}{\rr+2}}\big( 1+\frac{\mm}{\rr+1} \big).$
}
\end{m-lemma}\vskip-1ex
\begin{m-proof}
For the first inequality, $\uu_{\pp,\exc}$ is positive at $\xi_\blt=\big[(\rr+1)(\rr+2)\big]^{\frac{1}{\rr+2}}$ (see~\S\ref{sssct:down}). 
\end{m-proof}

For $\pp\approx1$, $x_0(\pp)$ is confined near $1$, so to get a lower envelope of $\vu{\pp,\exc}$, we may expand $\vu{\mm,+}^\pp$ about the origin. 
We have $\big[1+(\mm-1)s\big]^{-\frac{\pp}{\mm-1}}\les1-\pp s+\frac{\pp(\pp+\mm-1)s^2}{2}$, thus the function obtained by integrating 
$\vu{\mm,-}''=-x^\rr\Big[ 1-\frac{(2\mm-1)x^{\rr+2}}{(\rr+1)(\rr+2)}
+\frac{(2\mm-1)(3\mm-2) x^{2\rr+4}}{2(\rr+1)^2(\rr+2)^2}\Big],\,\vu{\mm,-}(0)=1,\,\vu{\mm,-}'(0)=0,$ 
is a lower envelope for $\uu_\exc^{(\pp)}$. The explicit formula is 
\\[.5ex]\centerline{
$\vu{\mm,-}(x)
=1-\frac{x^{\rr+2}}{(\rr+1)(\rr+2)}
+\frac{(2\mm-1)x^{2(\rr+2)}}{(\rr+1)(\rr+2)(2\rr+3)(2\rr+4)}
+\frac{(2\mm-1)(3\mm-2)x^{3(\rr+2)}}{2(\rr+1)^2(\rr+2)^2(3\rr+5)(3\rr+6)}.$
}\\[0ex] 

To estimate the precision of the envelopes, we return to  Lemma~\ref{lm:Uu}(iii--iv): the inequality there holds for $x<[(\rr+1)^2/2]^{\frac{1}{\rr+2}}$, which singles out $\Xi^-_\turn=\big[ \frac{\rr(\rr+2)}{\rr+2\pp} \big]^{\frac{1}{\rr+2}}$. A computation shows that $\aaa\big(\Xi^-_\turn(1)\big)=\big[ 1-\frac{\rr}{(\rr+1)^2} \big]^{\frac{1}{\rr+2}}$ and $\aph\big(\Xi^-_\turn(1)\big)=\frac{\rr}{\rr-1}$ are both about $1$. Therefore $\vu{\mm,\pm}$ is expected to estimate well $\uu_{\pp,\exc}$ on $[0,\Xi^-_\turn]$. 

Finally we construct a point $x_0^\mdl(\pp)\in[\xi_\blt,x_0^+]$ approximating $x_0(\pp)$, as we did in~\S\ref{ssct:x0}. (The factor $\aaa(\xi^+_\turn(1))=(\rr+1)^{-\frac{1}{\rr+2}}$ is included into $x_0^+$.)
\\[.5ex]\centerline{
$x_0^\mdl(\pp)=
0.6\cdot\big[(\rr+1)(\rr+2)\big]^{\frac{1}{\rr+2}}
+0.4\cdot\Big[ \frac{(\rr+1)^2(\rr+2)}{\rr+\mm+1} \Big]^{\frac{2}{\rr+2}} 
\cdot \Big[1+\frac{\mm}{(\rr+1)^{1+\frac{1}{\rr+2}}} \Big].
$
} 
\renewcommand{\arraystretch}{1.1}\textscale{.75}{
\begin{longtable}[c]{c|c||c|c|c|c|c|c|}
\hhline{~-------}
&\multicolumn{1}{c||}{\diagbox[width=6ex, height=5ex]{$\pp$}{$\rr$}}
&\textbf{1}
&\textbf{5}
&\textbf{10}
&\textbf{50}
&\textbf{100}
&\textbf{1000}
\\\hhline{~=======}
&\multirow{2}{*}{\textbf{0.5}}
&1.901$\,\mid\,$1.938
&1.745$\,\mid\,$1.759
&1.523$\,\mid\,$1.532
&1.1677$\,\mid\,$1.1699
&{\cellcolor{lgray}{1.0967$\,\mid\,$1.0978}}
&{\cellcolor{lgray}{1.01406$\,\mid\,$1.01418}}
\\
&
&1.14$\,\mid\,$1.6(-2)
&1.28$\,\mid\,$4.8(-3)
&1.22$\,\mid\,$1.7(-3)
&1.07$\,\mid\,$9.2(-5)
&1.04$\,\mid\,$2.4(-5)
&1.00$\,\mid\,$2.4(-7)
\\\hhline{~-------}
&\multirow{2}{*}{\textbf{0.9}}
&1.969$\,\mid\,$1.969
&1.778$\,\mid\,$1.773
&1.541$\,\mid\,$1.540
&1.1710$\,\mid\,$1.1715
&{\cellcolor{lgray}{1.0983$\,\mid\,$1.0987}}
&{\cellcolor{lgray}{1.01422$\,\mid\,$1.01426}}
\\
&
&1.02$\,\mid\,$8.8(-3)
&1.26$\,\mid\,$3.8(-3)
&1.21$\,\mid\,$1.4(-3)
&1.07$\,\mid\,$8.9(-5)
&1.04$\,\mid\,$2.3(-5)
&1.00$\,\mid\,$2.4(-7)
\\\hhline{--------}
\multicolumn{1}{|c|}{${\kern-.5ex}x_0(\pp)\mid x_0^\mdl(\pp){\kern-.5ex}$}
&\multirow{2}{*}{\textbf{1}}
&1.986$\,\mid\,$1.977
&1.787$\,\mid\,$1.776
&1.545$\,\mid\,$1.542
&1.1718$\,\mid\,$1.1719
&{\cellcolor{lgray}{1.0987$\,\mid\,$1.0989}}
&{\cellcolor{lgray}{1.01425$\,\mid\,$1.01428}}
\\ 
\multicolumn{1}{|c|}{$\Xi^-_\turn\mid\Delta\uu_\pm$}
&
&1.00$\,\mid\,$7.7(-3)
&1.25$\,\mid\,$3.6(-3)
&1.21$\,\mid\,$1.4(-3)
&1.07$\,\mid\,$8.8(-5)
&1.04$\,\mid\,$2.3(-5)
&1.00$\,\mid\,$2.4(-7)
\\\hhline{--------}
&\multirow{2}{*}{\textbf{1.01}}
&1.988$\,\mid\,$1.977
&1.788$\,\mid\,$1.777
&1.546$\,\mid\,$1.542
&1.1719$\,\mid\,$1.1720
&{\cellcolor{lgray}{1.0987$\,\mid\,$1.0989}}
&{\cellcolor{lgray}{1.01426$\,\mid\,$1.01429}}
\\ 
&
&0.99$\,\mid\,$7.6(-3)
&1.25$\,\mid\,$3.6(-3)
&1.21$\,\mid\,$1.4(-3)
&1.07$\,\mid\,$8.8(-5)
&1.04$\,\mid\,$2.3(-5)
&1.00$\,\mid\,$2.4(-7)
\\\hhline{~-------}
&\multirow{2}{*}{\textbf{1.1}}
&2.003$\,\mid\,$1.987
&1.795$\,\mid\,$1.780
&1.550$\,\mid\,$1.544
&1.1726$\,\mid\,$1.1723
&{\cellcolor{lgray}{1.0991$\,\mid\,$1.0991}}
&{\cellcolor{lgray}{1.01429$\,\mid\,$1.01430}}
\\ 
&
&0.97$\,\mid\,$6.9(-3)
&1.25$\,\mid\,$3.4(-3)
&1.20$\,\mid\,$1.4(-3)
&1.07$\,\mid\,$8.8(-5)
&1.04$\,\mid\,$2.3(-5)
&1.00$\,\mid\,$2.4(-7)
\\\hhline{~-------}
&\multirow{2}{*}{\textbf{1.5}}
&2.070$\,\mid\,$2.015
&1.829$\,\mid\,$1.794
&1.568$\,\mid\,$1.552
&1.1761$\,\mid\,$1.1740
&{\cellcolor{lgray}{1.1008$\,\mid\,$1.0999}}
&{\cellcolor{lgray}{1.01445$\,\mid\,$1.01438}}
\\ 
&
&0.90$\,\mid\,$4.7(-3)
&1.23$\,\mid\,$2.8(-3)
&1.20$\,\mid\,$1.2(-3)
&1.07$\,\mid\,$8.5(-5)
&1.04$\,\mid\,$2.3(-5)
&1.00$\,\mid\,$2.4(-7)
\\\hhline{~-------}
\caption{Error estimate $\Delta\uu_\pm=(\vu{\mm,+}-\vu{\mm,-})(\Xi^-_\turn)$ and value of first zero.}
\label{tab:x0smallp}
\end{longtable}
}\renewcommand{\arraystretch}{1}
\nit The shaded cells are boundary layers. The approximation $x_0(\pp)\!\approx\!x_0^\mdl(\pp)$ is precise for $\rr\!\ges\!10$.


\section{Applications and reflections}

We apply the results obtained so far. First, we verify the precision of our estimates for $x_0$ by solving~\eqref{eq:ef} backwards, starting from the approximate first zero, to see whether we recover the initial values at $x=0$. Second, we transform~\eqref{eq:ef} into a two-point boundary value problem, thus probing the validity of expressing $x_0$ in terms of $\pp,\rr$.

\subsection{Backward solution}\label{ssct:backw}

The knowledge of the approximate value of $x_0$ allows solving~\eqref{eq:ef} backwards. This implicitly means that we can estimate the value of $\uu_\exc'(x_0)$, too. We need one more point $\tld\xi_\mdl$, then consider the slope of the line joining $(x_0^\mdl,0)$ and $(\tld\xi_\mdl,\up(\tld\xi_\mdl))$: 
\\[0ex]\centerline{
$\; t:=\frac{0.44\cdot\pp+0.56\cdot\rr}{\pp+\rr+4}\in[0.4,0.6],\,
\tld\xi_\mdl:=(1-t)\cdot\xi^+_\turn+t\cdot\xi_\infl\in[\xi^+_\turn,\xi_\infl],\;\;
\mu:=\frac{-\up(\tld\xi_\mdl)}{x_0^\mdl-\tld\xi_\mdl}.$
}\\[0.5ex] 
(One can work out $\mu$ explicitly.) We define the backward-\eqref{eq:ef} equation as follows: 
\begin{align*}
\ww''=-x^\rr\ww^\pp,\;\;\ww(x_0^\mdl)=0,\;\ww'(x_0^\mdl)=\mu.\tag{EF*}
\label{eq:ef*}
\end{align*} 
Ideally, the output at $x=0$ should be $(1,0)$, the initial condition of~\eqref{eq:ef}. 
\setlength{\LTleft}{-2.5ex}\renewcommand{\arraystretch}{1.1}\textscale{.75}{
\begin{longtable}{c|c|c||c||c||c||c||c|}
\hhline{~-------} 
	&\multicolumn{1}{|c||}{\diagbox[width=6ex, height=5ex]{$\pp$}{$\rr$}}
    & \multicolumn{1}{c||}{\textbf{4}} & \multicolumn{1}{c||}{\textbf{15}} 
	& \multicolumn{1}{c||}{\textbf{60}} & \multicolumn{1}{c||}{\textbf{225}} 
	& \multicolumn{1}{c||}{\textbf{500}} & \multicolumn{1}{c|}{\textbf{1700}} 
\\ \hhline{-=======} 
\multicolumn{1}{|c|}{${\kern-.5ex}(\ww(0),\ww'(0)){\kern-.5ex}$}
	&\multicolumn{1}{|c||}{\multirow{2}{*}{\textbf{4}}}
	& (1.089, -0.043) & (0.990, 0.020)  & (0.875, 0.084) 
	& \cellcolor{lgray}{(0.898, 0.044)} & \cellcolor{lgray}{(1.066, -0.126)} & \cellcolor{lgray}{(1.970, -1.029)}
\\ 
\multicolumn{1}{|c|}{$\uu_\exc'(x_0)\mid\mu$}
&\multicolumn{1}{|c||}{}
& -1.1$\;\mid\;$-1.1 & -3.8$\;\mid\;$-3.9&-16.8$\;\mid\;$-16.3
& -67.0$\;\mid\;$-63.1 & -151.5$\;\mid\;$-141.1 & -521.9$\;\mid\;$-503.1   
\\ \cline{1-1}\hhline{~=======} 
	&\multicolumn{1}{|c||}{\multirow{2}{*}{{\textbf{15}}}} 
	& (1.033, -0.019) & (0.999, 0.019) & (0.922, 0.085) & (0.874, 0.130) 
	& \cellcolor{lgray}{(0.899, 0.106)}   
	& \cellcolor{lgray}{(1.177, -0.169)}
\\ 
&\multicolumn{1}{|c||}{}
& -0.4$\;\mid\;$-0.4 & -1.4$\;\mid\;$-1.5 & -6.2$\;\mid\;$-6.3
& -24.5$\;\mid\;$-24.8 & -55.3$\;\mid\;$-55.7 & -190.3$\;\mid\;$-193.8   
\\ \hhline{~=======} 
	&\multicolumn{1}{|c||}{\multirow{2}{*}{{\textbf{40}}}} 
	& (1.012, -0.011) & (1.004, -0.003) & (0.975, 0.028) 
	& (0.934, 0.058) & (0.932, 0.070) & \cellcolor{lgray}{(0.959, 0.043)}
\\ 
&\multicolumn{1}{|c||}{}
& -0.2$\;\mid\;$-0.2 & -0.6$\;\mid\;$-0.6 & -2.6$\;\mid\;$-2.6
& -10.1$\;\mid\;$-10.2 & -22.9$\;\mid\;$-23.1 & -78.7$\;\mid\;$-78.9
\\ \hhline{~=======} 
	&\multicolumn{1}{|c||}{\multirow{2}{*}{{\textbf{100}}}} 
	& (1.005, -0.006) & (1.004, -0.001) & (0.996, 0.005) 
	& (0.988, 0.012) & (0.987, 0.013) & (1.012, -0.011)
\\ 
&\multicolumn{1}{|c||}{}
& -0.1$\;\mid\;$-0.1 & -0.2$\;\mid\;$-0.2 & -1.1$\;\mid\;$-1.1
& -4.2$\;\mid\;$-4.2 & -9.5$\;\mid\;$-9.5 & -32.8$\;\mid\;$-32.8   
\\ \hhline{~=======} 
	&\multicolumn{1}{|c||}{\multirow{2}{*}{{\textbf{1000}}}} 
	& (1.0005, -0.0009) & (1.0006, -0.0004) & (1.0005, -0.0003) 
	& (1.0005, -0.0003) & (1.0005, -0.0004) & (1.001, -0.001)
\\ 
&\multicolumn{1}{|c||}{}
& -0.01$\;\mid\;$-0.01 & -0.03$\;\mid\;$-0.03 & -0.1$\;\mid\;$-0.1
& -0.4$\;\mid\;$-0.4 &-0.9$\;\mid\;$-0.9 & -3.3$\;\mid\;$-3.3
\\ \hhline{~-------} 
\caption{solving~\eqref{eq:ef} backwards.}
\label{tab:x0backw}
\end{longtable}
}\renewcommand{\arraystretch}{1}

\nit In the boundary layer cases, $\ww(0), \ww'(0)$ strongly depend on the initial values in~\eqref{eq:ef*}. The parameters $0.225$ in~\eqref{eq:manyx} and $0.56$ above make them small in these cases. (Further numerical tests show that the estimate $\uu'_\exc(x_0)\approx\mu$ is acceptable only for $\rr/\pp\les500$. The dependence on parameters becomes too strong.)

\subsection{Transforming the IVP into a BVP} \label{ssct:zx0}

The equation~\eqref{eq:ef} is an initial value problem, but one can change the viewpoint and transform it into an overdetermined boundary value problem. For $\zz>1$ given, one is interested to determine the (unknown) parameters $\rr,\pp$, such that the following equation admits solution: 
\begin{m-eqn}{
\biggl\{\begin{array}{l}
\uu''=-x^\rr\cdot \uu^\pp,\\ 
\uu(0)=1,\;\uu(\zz)=0,\;\uu'(0)=0.
\end{array}\biggr.
}\label{eq:overdet}
\end{m-eqn}
To our knowledge, there is no reference dealing with this matter. Yet it's easy justifying the interest in it: suppose one can measure the (physically relevant) first zero of $\uu$ but the parameters (e.g. inside a gas cloud) leading to it are unknown. Thus, the excessive condition $\uu'(0)=0$ is necessary due to physical considerations.

If it's omitted, usually there exists a unique solution (see~\cite{zhang,wei,xu,ch+ku,dai+li} for analytical approaches and \cite{wrd,sita} for numerical ones). In contrast, the vanishing of the derivative makes the problem overdetermined, it's solvable only for pairs $(\rr,\pp)$ satisfying some relation. This is the matter we wish to discuss. Obviously, numerical methods can not address this issue.

In the sequel, $\zz$ is the given desired zero in~\eqref{eq:overdet} and $x_0$ is the zero of $\uu_{\rr,\pp}$, where $\pp$ is a function of $\rr,\zz$. The smaller $|\zz-x_0|$ the better is the estimate of $\pp$ in terms of $\rr, \zz$. 

\begin{m-lemma}\label{lm:zx0} 
Let $(\rr,\pp)$ be a solution of~\eqref{eq:overdet}. The the following statements hold true: 
\begin{enumerate}[leftmargin=5ex]
\item 
We have the inequality $\rr>\rr_\mn$, where $\rr_\mn:=\min\{\rr\mid\zz>[(\rr+1)(\rr+2)]^{\frac{1}{\rr+2}}\}$.
\item The parameters $\rr, \pp$ are correlated as follows:
\begin{m-eqn}{
\pp\approx 2\biggl[\frac{\zz^{\frac{\rr+2}{\rr+1}}}{[(\rr+1)(\rr+2)]^\frac{1}{\rr+1}}-1\biggr]\cdot(\rr+1)^{1+\frac{1}{\rr+2}
}.\label{eq:prz}
}\end{m-eqn}
\end{enumerate}
\end{m-lemma}

It's unclear how to extract this kind of information from the literature. The difficulty is to have a \emph{sufficiently precise and solvable equation} relating  $\pp, \rr, \zz$. Thus the frequently used power/Puiseux series are not helpful. The simplicity of Newton's method is crucial.

\begin{m-proof}
(i) The estimate~\eqref{eq:u-blt} implies that $\zz>[(\rr+1)(\rr+2)]^{\frac{1}{\rr+2}},$ so necessarily $\rr>\rr_\mn$. 

\nit(ii) Since $x_0^\mdl$ is close to $x_0$, we use it to get an approximate relation between $\zz, \rr, \pp$: 
\\[.5ex]\centerline{\textscale{.9}{$
\zz
=\Big[ \frac{(\rr+1)(\rr+2)}{1+\frac{\mm}{\rr+1}} \Big]^{\frac{1}{\rr+2}}
\!\cdot 
\Big[ \lda+ \lda'\frac{\mm}{(\rr+1)^{\frac{\rr+3}{\rr+2}}}\Big],\;\lda,\lda'\approx1
\Leftrightarrow 
 \frac{1}{1+\frac{\mm}{\rr+1}}\cdot\Big[ \lda+ \lda'\frac{\mm}{(\rr+1)^{\frac{\rr+3}{\rr+2}}}\Big]^{\rr+2}
=\frac{\zz^{\rr+2}}{(\rr+1)(\rr+2)}.$}}\\[.5ex] 
At this point, we let $\lda=\lda'=1$ and replace $\frac{\mm}{\rr+1}\mt\frac{\mm}{(\rr+1)^{\frac{\rr+3}{\rr+2}}}$. The conclusion follows.
\end{m-proof}
Since $x_0^+>\zz$, the same argument shows that the following inequality is always true:   
\\[.5ex]\centerline{
$\pp> 2\biggl[\frac{\zz^{\frac{\rr+2}{\rr+1}}}{[(\rr+1)(\rr+2)]^\frac{1}{\rr+1}}-1\biggr]\cdot(\rr+1)-1.$
}\\[.5ex]
\nit The table below validates our estimate: the exact value $x_0$ is close to $\zz$, for each $(\rr,\pp)$. 
\renewcommand{\arraystretch}{1.1}\textscale{.75}{
\begin{longtable}[c]{|l|r||l|r||l|r||l|r|}
\cline{1-8} 
\multicolumn{2}{|c||}{$\zz=1.5,\;\;\rr_\mn=10$}&\multicolumn{2}{c||}{$\zz=1.1,\;\;\rr_\mn=94$}&
\multicolumn{2}{c||}{$\zz=1.01,\;\;\rr_\mn=1464$}&\multicolumn{2}{c|}{$\zz=1.001,\;\;\rr_\mn=19893$}  
\\ \hhline{========}
$\rr=15$&\multirow{2}{*}{$x_0=1.513$}&
$\rr=150$&\multirow{2}{*}{$x_0=1.097$}&
$\rr=1600$&\multirow{2}{*}{$x_0=1.0099$}&
$\rr=20500$&\multirow{2}{*}{$x_0=1.00100$}
\\ 
$\pp=4$&&$\pp=10$&&
$pp=3$&&$\pp=2$&
\\ \hhline{========}
$\rr=150$&\multirow{2}{*}{$x_0=1.492$}&
$\rr=2500$&\multirow{2}{*}{$x_0=1.098$}&
$\rr=7000$&\multirow{2}{*}{$x_0=1.0096$}&
$\rr=25000$&\multirow{2}{*}{$x_0=1.00098$}
\\ 
$\pp=128$&&$\pp=468$&&
$\pp=105$&&$\pp=10$&
\\  \hhline{========} 
$\rr=1500$&\multirow{2}{*}{$x_0=1.497$}&
$\rr=7000$&\multirow{2}{*}{$x_0=1.099$}&
$\rr=35000$&\multirow{2}{*}{$x_0=1.0098$}&
$\rr=70000$&\multirow{2}{*}{$x_0=1.00097$}
\\ 
$\pp=1466$&&$\pp=1364$&&
$\pp=658$&&$\pp=96$&
\\  \hhline{========}
$\rr=10000$&\multirow{2}{*}{$x_0=1.499$}&
$\rr=25000$&\multirow{2}{*}{$x_0=1.099$}&
$\rr=70000$&\multirow{2}{*}{$x_0=1.0099$}&
$\rr=150000$&\multirow{2}{*}{$x_0=1.00098$}
\\ 
$\pp=9956$&&$\pp=4958$&&
$\pp=1355$&&$\pp=253$&
\\ \cline{1-8}
\caption{Testing~\eqref{eq:prz} for $\zz-1$ small.}\label{tab:z-1small}
\end{longtable}
\begin{longtable}[c]{|l|r||l|r||l|r||l|r|}
\cline{1-8} 
\multicolumn{2}{|c||}{$\zz=3$}&
\multicolumn{2}{c||}{$\zz=50$}&
\multicolumn{2}{c||}{$\zz=250$}&
\multicolumn{2}{c|}{$\zz=500$}  
\\  \hhline{========}
$\rr=10$&\multirow{2}{*}{$x_0=3.03$}&
$\rr=10$&\multirow{2}{*}{$x_0=51.48$}&
$\rr=10$&\multirow{2}{*}{$x_0=258.00$}&
$\rr=10$&\multirow{2}{*}{$x_0=516.08$}
\\ 
$\pp=31$&&$\pp=1202$&&
$\pp=7091$&&$\pp=15135$&
\\  \hhline{========} 
$\rr=50$&\multirow{2}{*}{$x_0=2.99$}&
$\rr=50$&\multirow{2}{*}{$x_0=50.72$}&
$\rr=50$&\multirow{2}{*}{$x_0=253.87$}&
$\rr=50$&\multirow{2}{*}{$x_0=507.79$}
\\ 
$\pp=179$&&$\pp=4978$&&
$\pp=26148$&&$\pp=53125$&
\\  \hhline{========} 
$\rr=300$&\multirow{2}{*}{$x_0=2.99$}&
$\rr=300$&\multirow{2}{*}{$x_0=50.14$}&
$\rr=300$&\multirow{2}{*}{$x_0=250.78$}&
$\rr=300$&\multirow{2}{*}{$x_0=501.58$}
\\ 
$\pp=1165$&&$\pp=29305$&&
$\pp=149782$&&$\pp=300872$&
\\  \hhline{========} 
$\rr=1000$&\multirow{2}{*}{$x_0=2.99$}&
$\rr=1000$&\multirow{2}{*}{$x_0=50.04$}&
$\rr=1000$&\multirow{2}{*}{$x_0=250.23$}&
$\rr=1000$&\multirow{2}{*}{$x_0=500.48$}
\\ 
$\pp=3955$&&$\pp=97784$&&
$\pp=497787$&&$\pp=998282$&
\\ \cline{1-8} 
\caption{Testing~\eqref{eq:prz} for $\zz$ large.}
\label{tab:z-1large}
\end{longtable}
}\renewcommand{\arraystretch}{1}


\subsection{Reflections, speculations}\label{ssct:spec}

\subsubsection{} The naively defined slope $\mu$ approximates well $\uu_\exc'(x_0)$, except the `strong boundary' cases (top-right in Table~\ref{tab:x0backw}). Roughly, it means that the knowledge of $x_0$ determines the value of the derivative. This is somewhat surprising, especially when trying to interpret in physical terms. The function $\uu_\exc$ represents (after coordinate changes) the density of a star, so $x_0$ is its radius (density becomes zero). The derivative $\uu'(x_0)$ is the rate of change of the density at the boundary of the star. The estimate $\uu'(x_0)\approx\mu$ means that (gaseous) stars/polytropes are able to `guess' their own matter's vanishing rate at the boundary, as soon as they know their size. Should this argument be correct, a physical explanation would be welcome.

\subsubsection{} We reached the \eqref{eq:ef}-equation (see~\S\ref{sssct:ode}) starting from $\yy''=\kp^2\yy^\nn,\,\yy(0)=1$, with relevance to radiative heat transfer (the classical Stefan-Bolzmann law corresponds to $\nn=4$). The approximate solutions of these ODEs are related by a simple change of variable. 
Since the equations describe physical realities, one naturally wonders if there is any physically motivated connection between self-gravitating gaseous spheres and heat radiation. 
After examining the literature, the author ---certainly not a physicist--- suspects that the link is given by Tsallis' non-extensive statistics leading in both situations to power-law formulae. 
\begin{itemize}[leftmargin=3ex] 
\item 
Tsallis' entropy adequately describes~\cite{plas} polytropic models. The entropy is maximized by self‑gravi\-ta\-ting systems, the resulting equilibrium states are stellar polytropes~\cite{sak+tar,tar+sak}.

\item 
Tsallis statistics is used~\cite{len+men} to deduce a multi-dimensional generalization of the blackbody radiation formula which includes non-extensive systems possessing a large number of degrees of freedom. The exponent (power) appearing in the generalized, non-extensive Stefan-Bolzmann law~\cite{mppt} is not necessarily four any more [\textit{id.}, Fig.\,2].
\end{itemize}


\section{Conclusion}

We investigated the dependence on the parameters $\rr,\pp$ of the exact solution $\uu_\exc$ of the following well-known Emden-Fowler equation: 
$$
\begin{array}{rcl}
{\rm(EF)}:&\uu''=-x^\rr\cdot\uu^\pp,& \uu(0)=1,\; \uu'(0)=0.
\end{array} 
$$ 
To achieve our goal, we proceeded as follows: 
\begin{enumerate}[leftmargin=5ex]
\item 
We set up a basic analytical framework which allows constructing upper and lower envelopes. These are used to deduce qualitative information about: 
\begin{itemize}[leftmargin=3ex]
\item 
approximate values of the exact solution and its derivative;
\item 
conditions leading to boundary layer for the solution; 
\item 
the location of the boundary-layer-interval. 
\end{itemize}
\item 
We apply the analytical tools to investigate the geometric properties of the exact solution $\uu_\exc$ of (EF). A central role in our approach is played by the turning point $\xi_\turn$, where the graph of $\uu_\exc$ bends downwards the most; it satisfies the equation $\uu_\exc'''(x)=0$. We determine upper/lower envelopes of $\uu_\exc$, which allow estimating $\xi_\turn$ and the slope $\uu_\exc'(\xi_\turn)$. Using them, we apply Newton's method to deduce bounds for the first zero of $\uu_\exc$, in terms of $\rr,\pp$. We verify our results by backward solving~\eqref{eq:ef}. Furthermore, given $\zz>1$ ---the value of the desired first zero---, we approximate $\pp$ as function of $\rr, \zz$, so that $\uu_\exc$ vanishes at $\zz$. 
\end{enumerate} 


\end{document}